\documentclass[11pt,a4paper]{article}
\usepackage{epsfig}
\usepackage{graphicx}
\usepackage{amsmath}
\usepackage{amssymb}

\newcommand{\bbN}{{\mathbb N}}

\title{Modelling the spatial organization of cell proliferation
in the developing central nervous system}
\author{
Jean Clairambault \thanks{INRIA-Rocquencourt, projet BANG, Domaine de Voluceau, BP105, F 78153
LeChesnay cedex, \texttt{Jean.Clairambault@inria.fr}}
\and
Vladimir Flores \thanks{Interdisciplinary Group in Theoretical Biology, Department of Biostructural
Sciences, Favaloro University. Sol\'{i}s 453 (1078),
Buenos Aires, Argentina, \texttt{vflores@favaloro.edu.ar}}
\and 
Beno\^{i}t Perthame\footnotemark[1] \thanks{\texttt{benoit.perthame@upmc.fr}}
\and 
Melina Rapacioli\footnotemark[2] \thanks{\texttt{mrapacioli@favaloro.edu.ar}}
\and
Edmundo Rofman  \thanks{IAM-CONICET, Argentine and INRIA, France, \texttt{edmundo.rofman@wanadoo.fr}}
\and
Rafael Verdes   \thanks{Instituto de Matematica ''Beppo Levi'', Universidad Nacional de Rosario,
Argentina, \texttt{ rverdes@arnet.com.ar}}
}

\begin{document}
\maketitle

\abstract{How far is neuroepithelial cell proliferation in the developing central nervous system a deterministic process? Or, to put it in a more precise way, how accurately can it be described by a deterministic mathematical model? To provide tracks to answer this question, a deterministic system of transport and diffusion partial differential equations, both physiologically and spatially structured, is introduced as a model to describe the spatially organized process of cell proliferation during the development of the central nervous system. As an initial step towards dealing with the three-dimensional case, a unidimensional version of the model is presented. Numerical analysis and numerical tests are performed. In this work we also achieve a first experimental validation of the proposed model, by using cell proliferation data recorded from histological sections obtained during the development of the optic tectum in the chick embryo.}

\section{Introduction}
The present work aims at modeling the spatial organization of the neuroepithelial (NE) cell proliferation in a developing cortical structure along an early and brief developmental period by using sets of quantitative data empirically obtained from a standardized experimental model: the developing chick optic tectum (OT). The mathematical model is based on a deterministic approach that uses the formalism of partial differential equations.\\

\subsection{Biological background: {\bf Developmental neurobiology}}

\subsubsection{\bf Relevance of a spatial and temporal organization in a developing system}

The appropriate number of cells of each terminally differentiated cell type and also the spatial patterns they exhibit within the different tissues and organs composing pluricellular organisms are governed by interactive self-regulating behaviors that the developing cells exhibit during the embryonic development \cite{Patternformation}. The increase in supra-cellular complexity generated during development requires the temporally and spatially organized operation of specific developmental cell behaviors (DCBs). These DCBs are usually reciprocally regulated and operate simultaneously and interactively \cite{Gilbert}. Every developing cell population can be considered as both emitter and receiver of developmental regulatory signals having, on the one hand, informative and, on the other hand, structural roles. Thus, a central hypothesis in Developmental Biology proposes that the space-time organized operation of developmental cell behaviors depends on the cooperative establishment of spatially organized cell signaling networks mediated by diffusing informative molecules \cite{[1], [2], [3], [4]}. Molecular diffusion results in asymmetric distribution of developmentally active informative signals. This asymmetry plays a fundamental role in establishing temporal and/or spatial organization of specific DCBs that result in the whole developmental process in the organized patterns of cells, tissues and organs, that living organisms eventually exhibit in their final, terminally differentiated, state.\\

\subsubsection{\bf The organizers and CNS patterning}

Not all developing cell populations possess similarly relevant informative roles. There exist specific transient cell populations, the so-called organizers, that primarily play informative roles influencing or regulating the developmental behavior of the other cells. By means of installing asymmetric distributions of developmentally active signals the organizers serve as instruments of an informative reference system in the establishment of spatially organized processes of cell determination and differentiation \cite{[Schoenwolf], [Towers]}. Amongst the most complex biological structures, the multilayered concentric neuronal organization of the central nervous system (CNS), i.e., brain cortex, cerebellum cortex etc., occupy a privileged position. The development of such a structural and functional complexity, the so-called corticogenesis, requires the organized operation of several DCBs. Amongst these DCBs: (a) the cell proliferation (CP), counteracted by apoptosis or programmed cell death, is involved in the generation of the appropriate number of neurons for each cortical area and each cortical layer; (b) the directed cell migration, a process mediated by specific interfacial interactions between cell surface and extracellular matrix components, controls the correct position of each specific neuronal type along the CNS spatial axes; (c) cell determination or commitment, a process mediated by irreversible genetic information reprogramming, allows different cell populations to select one out of a set of multiple developmental pathways; (d) cell differentiation, a process mediated by selective gene activation and selective protein synthesis, warrants the expression of different specific neuronal phenotypes; finally, (e) the development of the neural processes (neuritogenesis), dendrites and axons and the establishment of specific synaptic contacts (synaptogenesis) conclude this complex self-organizing and interactively regulated process of corticogenesis.\\

All these DCBs are properly organized in time and space. During the early
stages, several organizers strategically positioned along the cephalic-caudal and dorsal-ventral axes of the developing CNS establish the primitive pattern of each region and subregion, as well as the identity of the neuromeres. Later, cell proliferation, migration and differentiation controlled by reciprocal interactions allow the whole system to expand and differentiate into several structurally and functionally integrated systems of interconnected neuronal circuits. Abundant evidence indicate that all these processes are spatially and temporally organized \cite{[Brusco95], [Scicolone95], [PereyraAlfonso97], [PereyraAlfonso98], [Sanchez02], [RodGil05]}.\\

\subsubsection{\bf Roles of neuroepithelial cell proliferation and neuronal migration in corticogenesis}

During the early development, the primitive CNS primordium, the so-called
neural tube (NT), is almost exclusively composed of proliferative
neuroepithelial (NE) cells. During this early proliferative phase the NE
cells behave as a population of self-renewing stem cells that divide symmetrically, i.e. from each dividing NE cell originates two similar NE
cells. These dividing NE cells are typically located along the inner surface of the NT forming the so-called generation zone (GZ). This type of early proliferation produces an amplification of the NE cells: they significantly increase in number, and the GZ increases in both area (planar expansion) and thickness (radial expansion). During a later proliferative phase NE cells proliferate asymmetrically, i.e. from each dividing NE cell originates a NE cell and a postmitotic neuron (PN). The number of NE cells and the GZ size stabilize while postmitotic neurons accumulate at an underlying premigratory zone (PMZ). Later, neurons migrate radially and the PMZ weakens while different layers of concentrically organized neurons begin to organize the future cortex. By the end of the proliferative phase the number of NE cells decreases, the GZ becomes thinner and finally disappears.\\

\subsubsection{\bf Neuroepithelial cell proliferation analyzed as a stochastic point process}

The NE cells compose a synchronized and functionally integrated population of stem cells. Biomathematical analyses of numerical series representing the position of mitotic NE cells along the OT cephalic-caudal axis, analyzed within the framework provided by the theory of stochastic point processes, show that NE cells do not proliferate independently. It can be shown that the proliferative dynamics embeds (subsumes) two different components: (a) a non-stationary component that can be deterministically described and reveals a spatially organized process and (b) a stationary stochastic point process representing a uniformly distributed basal proliferative activity \cite{[Rapacioli00], [Mazzeo04a]}. More recently a different approach allows demonstrating that the stochastic component corresponds to an anti-correlated stationary process implying the result that there are short-range inhibitory interactions between neighboring proliferating NE cells \cite{[Mazzeo08], [Rapacioli08]}. These mathematical results are consistent with experimental evidence indicating that cell proliferation is a DCB submitted to cell-cell controlling interactions and also to intracellular controlling processes. 

\subsubsection{\bf Extra- and intracellular signals controlling proliferation dynamics}

A general notion establishes that the coordinated behavior shown by proliferating cells during development \cite{segmentation} depends, on one hand, on a set of extracellular signals, the so-called growth factors \cite{[Raff], [Bohnsack]} and a set of intracellular informative signals composing a cell cycle control center whose main members are the regulatory proteins cyclins as well as cyclin-dependent kinases \cite{Johnson}.\\

The temporally organized operation of these regulatory proteins installs the so-called cell division cycle (proliferative cycle). The cell division cycle is composed of a periodic ordered sequence of cyclic intracellular molecular events with typical durations. It involves two main phases: (a) the mitotic or M phase that includes the karyokinesis (nuclear division) and the cytokinesis (cytoplasm division) and (b) the interphase, or intermitotic phase. In turn, each of these phases involves temporally ordered sub-phases. The M phase is composed of the pro-, prometa-, meta-, ana- and telo-phase while the interphase is composed of the G1, S and G2 phases. Each type of proliferative cells population displays typical cell cycle length oscillating around a mean value. The total cell cycle length and also those of the different phases, vary, depending on the cell type. Besides, particular cell types may change their cell cycle durations when long time series of subsequent periods are considered.\\

\subsection{Biological settings of the present study}

\subsubsection{\bf The experimental model system: NE cell proliferation in the developing avian optic tectum}

The avian optic tectum (OT) is a mesencephalic alar plate derivative. The main OT input proceeds from the retinal ganglion cell axons. This information is integrated into others inputs, processed and transferred to other CNS areas.
The OT possesses a typical cortical, multilayered concentric organization, consisting of alternating neuronal and fibrous layers. The developing avian OT is a well-characterized experimental model of CNS development. It is known that the OT determination and the positional information specifying the OT patterning are submitted to the organizing influence of the Isthmic Organizer (IsO). The IsO organizing activity is mediated by informative molecules which specify positional information along the OT cephalic-caudal axis [9, 10, 11, 12]. The positional information is firstly established as asymmetric gene expression and is then ``translated'', by the differential operation of several DCBs into a gradient of cytoarchitectonic differentiation. There are well-documented data about its proliferative kinetics \cite{[LaVail], [Rapacioli00], [Mazzeo04a], [Rapacioli08]}, the postmitotic neuronal migration \cite{[Rapacioli00]} and the expression of enzymes involved in neuronal migration \cite{[PereyraAlfonso97]}, the histogenesis of its multilaminated architectural pattern \cite{[Scicolone95]} and also on the temporal and spatial developmental pattern of its specific neuronal types \cite{[Sanchez02], [Scicolone06]}. There are also detailed structural, pharmacological and molecular studies on the development of its extrinsic innervation (serotonin innervation, retino-tectal connections) \cite{[Brusco95], [Shigetani], [Luksch]}. The developmental patterns of expression of several molecular components involved in synaptogenesis of the OT local circuits are also well characterized \cite{[RodGil05], [Vacotto], [Fiszer]}.\\

\subsubsection{\bf Biological methods and cell proliferation records}

Empirical data on NE cell proliferation analyzed in this paper were recorded from complete histological serial sections obtained from OT of different developmental stages. OTs obtained from 2, 4 and 6 embryonic (E) days ($E2$, $E4$, $E6$) were used. During this interval the OT undergoes significant changes in size and shape. After dissection, OTs were processed for conventional histological methods \cite{[Scicolone95]}. Specimens were spatially oriented in order to obtain planar sections fulfilling criteria of adequacy between planar section orientation and developmental gradient axis position. These criteria were given in [6]. Cell proliferation records (CPR) along the OT cephalic-caudal axis were constructed on a video camera screen (Axioplan 2 imaging optical epifluorescence microscope coupled to an Axiocam HR color digital scanner and a computer equipped with the software Axiovision (all Carl Zeiss, Germany)) at a final magnification of 1000 X. A CPR is defined as a numerical sequence that indicates the density of mitotic cells in successive ordered, $25\mu m$ length windows, along OT cephalic-caudal axis. The total number of cells in the GZ was also computed for each $25 \mu m$ window.\\

\section{Complete mathematical model}

The present modeling concerns the development of the OT between embryonic days $E2$ and $E6$. It only describes the elongation of the cephalic-caudal axis and thus is only one space dimension. It also assumes that:
\\
(a) the NE cells divide symmetrically (the number of asymmetric division is negligible), 
\\
(b) the number of postmitotic neurons in the GZ is negligible and the increase in cell number in the GZ is mainly due to the symmetric divisions and 
\\
(c) the cell cycle length does not change significantly neither along the interval $E2$-$E4$ nor between $E4$ and $E6$. Consequently the numerical analysis and experimental validation treat separately these two periods. 
\\

We follow a classical description of the cell cycle that can be found in \cite{[B3]} and references therein. Our model is based on equations for the cell population densities structured by their progression along each phase $0 \leq a \leq 1$ and their position $x$ along the cephalic-caudal axis. The cells are divided in subclasses, as follows:
\\
\\
$n_{1}(t,a,x)=$linear density of neuroepithelial cells in phases $G_{1},S,G_{2}$,
\\[2mm]
$n_{2}(t,a,x)=$linear density of neuroepithelial cells in phase $M$,
\\[2mm]
$n_{3}(t,a,x)=$ postmitotic neurons generated by asymmetric mitosis. These
cells stop proliferating; they are destined to migrate and later differentiate
into specialized neurons.
\\

We propose a unique model for time $t$ between days $E2$ and $E10$. Cells in the phases are also subject to random motion along the cephalic-caudal axis and this leads us to represent the dynamics of cell division and motion by the equations \begin{equation} \label{(1)}
\left\{
\begin{array}{l}
\frac{\partial n_{1}}{\partial t}+\frac{\partial \lbrack v_{1}(t)n_{1}]}{
\partial a}=\frac{\partial ^{2}[D(x)n_{1}]}{\partial x^{2}},
\\[2mm]
\frac{\partial n_{2}}{\partial t}+\frac{\partial \lbrack v_{2}(t)n_{2}]}{
\partial a}=\frac{\partial ^{2}[D(x)n_{2}]}{\partial x^{2}},
\\[2mm]
\frac{\partial n_{3}}{\partial t}=q(t,x)v_{2}(t)n_{2}(t,a=1,x) .
\end{array}
\right.
\end{equation}

This system is completed with boundary and initial conditions
\begin{equation} \label{(1bc)}
\left\{
\begin{array}{l}
v_{2}(t)n_{2}(t,a=0,x)=v_{1}(t)n_{1}(t,a=1,x), 
\\[2mm]
v_{1}(t)n_{1}(t,a=0,x)=[2-q(t,x)]v_{2}(t)n_{2}(t,a=1,x), 
\\[2mm]
\frac{\partial n_{1}}{\partial x}(t,a,x=x_{+}(t))=0, \quad \frac{\partial n_{2}}{
\partial x}(t,a,x=x_{+}(t))=0, \quad  \frac{\partial n_{3}}{\partial x}
(t,a,x=x_{+}(t))=0, 
\\[2mm]
n_{1}(0,a,x)=n_{1}^{0}(a,x), \qquad n_{2}(0,a,x)=n_{2}^{0}(a,x), \qquad 
n_{3}(0,a,x)=0.
\end{array}
\right.
\end{equation}

We have used the following notations 
\\
\\
$\bullet$ $t$ the chronological time variable, 
\\[2mm]
$\bullet$ $x$: abscissa of transverse sections along the cephalic-caudal OT axis, $0\leq x\leq x_{+}(t)$,  $x=0$ cephalic tip, $x=x_{+}(t)$, caudal tip at time $t$,
\\[2mm]
$\bullet$ $0\leq a\leq 1$: degree of progression along the mitotic ($M$) phase or the intermitotic 
$(G_{1}, S, G_{2})$ phase,
\\[2mm]
$\bullet$ $v_{i}(t),\;(i=1,2)$, homogeneous to an inverse of time, is the frequency of the corresponding phase duration,
\\[2mm]
$\bullet$ $D(x)$: ``effective cellular diffusion'' at section $x$,
\\[2mm]
$\bullet$ $0\leq q(t,x)\leq 1$: postmitotic neuronal population growth rate, or rate of differentiation (a more complete model should consider it as a nonlinear function of the cell densities).
\\

The model should also be completed with a law for the OT elongation along the development, which means a description of the 
$x_{+}(t)$. It should naturally take the form
\begin{equation} \label{eqmotion}
\frac{d}{dt} x_+(t)=  F\left(n_1(t), n_2(t), n_3(t) \right) ,
\end{equation}
with a function $F(\cdot)$ that describes how the OT development is ruled by the different cell types. In this case the  Neumann boundary conditions at $x^+(t)$ should be replaced by the no-flux conditions 
$$
\dot x_+(t) n_i(t,a,x_+(t)) + D(x) \frac{\partial}{\partial x} n_i(t,a, x=x_+(t))=0.
$$

To access such nonlinear coefficients as $q(\cdot)$ and $F(\cdot)$ is too demanding in a first stage, and thus simplifications are needed. Therefore we reduce the model complexity with two assumptions. Firstly, for days between $E2$ and $E6$ we can neglect asymmetric cell proliferation which leads to $q=0$ and $n_3=0$. Secondly, we assume that the OT has constant length between days $E2$-$E4$ and $E4$-$E6$; this leads to fix $x_+$ independent of time with a jump at the end of day $E4$.

\section{Reduced model for the symmetric proliferative phase
without postmitotic migrating neurons ($q = 0$)}

We reduce the model complexity in the present work with two assumptions. Firstly, for days between $E2$ and $E6$ we can neglect asymmetric cell proliferation which leads to $q=0$ and $n_3=0$. Secondly, we assume that the OT has constant length between days $E2$-$E4$ and $E4$-$E6$; this leads to fix $x_+$ independent of time with a jump at the end of day $E4$. According to these assumptions, we consider only a simplified version of the model that we present now. We also present its discretization.

\subsection{The simplified model}

We define the region
\[
\Omega := \big\{ 0\leq t\leq T,\quad  0\leq a\leq 1,\quad  0\leq x\leq x_{+} \big\}, 
\]
with $x_+$ constant and we also suppose that $v_{1},v_{2},D$ constant. This leads to  
the model defined by the system of partial differential equations
\begin{equation}
\label{(2)}
\left\{
\begin{array}{l}
\frac{\partial n_{1}}{\partial t}+v_{1}\frac{\partial n_{1}}{\partial a}=D
\frac{\partial ^{2}n_{1}}{\partial x^{2}} ,
\\[2mm]
\frac{\partial n_{2}}{\partial t}+v_{2}\frac{\partial n_{2}}{\partial a}=D
\frac{\partial ^{2}n_{2}}{\partial x^{2}} ,
\end{array}
\right.
\end{equation}
completed with the boundary and initial conditions 
\begin{equation}
\label{(2bc)}
\left\{
\begin{array}{l}
v_{2}n_{2}(t,a=0,x)=v_{1}n_{1}(t,a=1,x), 
\\[2mm]
v_{1}n_{1}(t,a=0,x)=2v_{2}n_{2}(t,a=1,x), 
\\[2mm]
\frac{\partial n_{1}}{\partial x}(t,a,x=x_{+}(t))=0, \quad \frac{\partial n_{2}}{
\partial x}(t,a,x=x_{+}(t))=0, 
\\[2mm]
n_{1}(0,a,x)=n_{1}^{0}(a,x), \qquad n_{2}(0,a,x)=n_{2}^{0}(a,x).
\end{array}
\right.
\end{equation}

For later purposes we impose the compatibility condition between the boundary and initial conditions:
\begin{equation}\label{(3)}
v_{1}n_{1}^{0}(1,x)=v_{2}n_{2}^{0}(0,x)
\end{equation}
It implies that the solution is smooth enough for our analysis below, see \cite{[3]}.

\subsection{Discretization by the method of finite differences}

For every $K\in \bbN, \,I\in \bbN,\,J\in\bbN$ we define a mesh in $\Omega$ by the points
\[
(k\Delta x,i\Delta a,j\Delta t), \ \ k=0,1,...,K , \ i=0,1,...,I,\  j=0,1,...,J
\]
with 
$$
\Delta t=\frac{T}{K}, \quad  \Delta a=\frac{1}{I} , \quad \Delta x=\frac{x_{+}}{J}.
$$
For explicit schemes, these parameters have to be chosen with the stability (Courant-Friedrichs-Lewy) condition
\begin{equation}\label{14}
v\frac{\Delta t}{\Delta a}+2D\frac{\Delta t}{(\Delta x)^{2}}\leq 1.
\end{equation}

According to well established methods \cite{[B1],[B2]}, we approximate the derivatives in the differential equations by upwind (for first order derivatives) or centered (for the second order derivatives) differences on the mesh points. For a generic point of the mesh $(k\Delta x,i\Delta
a,j\Delta t)$ we set as usual 
\[
\frac{\partial n_{l}}{\partial t}\cong \frac{
n_{l}(k+1,i,j)-n_{l}(k,i,j)}{\Delta t},\ \ \ \ \ \ \ \frac{\partial
n_{l}}{\partial a}\cong \frac{n_{l}(k,i,j)-n_{l}(k,i-1,j)}{\Delta a}
\]

\[
\frac{\partial ^{2}n_{l}
}{\partial x^{2}}\cong \frac{n_{l}(k,i,j+1)-2n_{l}(k,i,j)+n_{l}(k,i,j-1)}{%
(\Delta x)^{2}}, \ \ \ \ l=1,2.
\]

In the following we will denote by:
\\
$\bullet$ $n(k,i,j)$ \ value of the exact solution at the point \textit{(k}$\Delta
x,i\Delta a,j\Delta t)$ of the mesh,
\\
$\bullet$ $u(k,i,j)$ \ \ value of the discretized solution at the point $(k)\Delta x,i\Delta a,j\Delta t)$,
\\
and we consider the system of equations obtained upon application of these approximations, with the corresponding boundary and initial conditions.
\\

The discretized model reads, for $i=1,...,I$, $j=1,...,J-1$, $k=0,1,...,K$, 
\begin{equation} \label{(4)}
\left\{
\begin{array}{l}
\frac{u_{1}(k+1,i,j)-u_{1}(k,i,j)}{\Delta t}+v_{1}\frac{
u_{1}(k,i,j)-u_{1}(k,i-1,j)}{\Delta a}=D\frac{
u_{1}(k,i,j+1)-2u_{1}(k,i,j)+u_{1}(k,i,j-1)}{(\Delta x)^{2}} ,
\\[2mm]
\frac{u_{2}(k+1,i,j)-u_{2}(k,i,j)}{\Delta t}+v_{2}\frac{
u_{2}(k,i,j)-u_{2}(k,i-1,j)}{\Delta a}=D\frac{
u_{2}(k,i,j+1)-2u_{2}(k,i,j)+u_{2}(k,i,j-1)}{(\Delta x)^{2}} ,
\\[2mm]
v_{2}u_{2}(k,0,j)= v_{1}u_{1}(k,I,j),  
\\[2mm]
v_{1}u_{1}(k,0,j)=2v_{2}u_{2}(k,I,j) ,
\\[2mm]
u_{1}(0,i,j)=n_{1}^{0}(i,j),\qquad u_{2}(0,i,j)=n_{2}^{0}(i,j).
\end{array}
\right.
\end{equation}

It is again convenient to also assume a compatibility condition between the boundary and initial conditions in the discretized model
\begin{equation}\label{(5)}
v_{1}n_{1}^{0}(I,j)=v_{2}n_{2}^{0}(0,j) \ \ \ \ j=0,1,...,J.
\end{equation}

We use $u$ to indicate both $u_{1}$ and $u_{2}$,
and $v$ for $v_{1}$ and $v_{2}$. Therefore, from equation \eqref{(4)}, we deduce 
\[
\begin{array}{ll}
\dfrac{u(k+1,i,j)-u(k,i,j)}{\Delta t}&+\,v\,\dfrac{u(k,i,j)-u(k,i-1,j)}{\Delta a}=
\\[2mm]
&= D\dfrac{u(k,i,j+1)-2u(k,i,j)+u(k,i,j-1)}{(\Delta x)^{2}}.
\end{array}
\]
In other words, we obtain $u(k+1,i,j)$ thanks to the iterations on the time label
\begin{equation}
\label{(6)}
\begin{array}{ll}
u(k+1,i,j)=& u(k,i,j)-v\dfrac{\Delta t}{\Delta a}(u(k,i,j)-u(k,i-1,j))+
\\[2mm]
&+D\dfrac{\Delta t}{(\Delta x)^{2}}(u(k,i,j+1)-2u(k,i,j)+u(k,i,j-1)).
\end{array}
\end{equation}
We use the discrete Neumann boundary conditions  
\begin{equation}\label{(7)}
u(k,i,0)=u(k,i,1) , \qquad  u(k,i,J+1)=u(k,i,J),
\end{equation}
for $k=0,1,...,K, \ \ i=0,1,...,I$.

Applying the induction formula \eqref{(6)} to calculate the values of $u_{1}$
and $u_{2}$, and bearing \eqref{(7)} and the initial conditions in mind, we obtain the values of the discretized solution
\begin{equation}
\label{(8)}
u_{1}(k,i,j),\ \ u_{2}(k,i,j), \ \ k=0,1,...,K, \ \ i=0,1,...,I, \ \ j=0,1,...,J.
\end{equation}

This describes the numerical method we use in the sequel together with experimental data.

\subsection{Convergence of the approximate solutions to the exact solution on the points of the mesh}

The convergence analysis of this numerical scheme is standard and can be found in \cite{[B1],[B2]}. We give here a fast account for the sake of completeness.

We use the notation $n(k,i,j)$ to indicate both $n_{1}(k,i,j)$ and $n_{2}(k,i,j)$,
and also $u(k,i,j)$ for $u_{1}(k,i,j)$ and $u_{2}(k,i,j)$.

Applying Taylor's formula to the exact solution at the points of the mesh we obtain:

\begin{equation}
\label{9}
\begin{array}{ll}
\dfrac{n(k+1,i,j)-n(k,i,j)}{\Delta t}&+\,v\,\dfrac{n(k,i,j)-n(k,i-1,j)}{\Delta a}=
\\[2mm]
&=D\dfrac{n(k,i,j+1)-2\,n(k,i,j)+n(k,i,j-1)}{(\Delta x)^{2}}+\Psi,
\end{array}
\end{equation}
where for some constants $A,B,C$, we have 
\begin{equation}
\label{10}
\left| \Psi \right| \leq A\Delta t+B\Delta a+C(\Delta x)^{2}.
\end{equation}

From \eqref{9} we obtain
\begin{equation} \label{11}
\begin{array}{ll}
n(k+1,i,j)=&n(k,i,j)-v\frac{\Delta t}{\Delta a}%
(n(k,i,j)-n(k,i-1,j))+\vspace{2mm}\\
&+D\frac{\Delta t}{(\Delta x)^{2}}(n(k,i,j+1)-2%
\,n(k,i,j)+n(k,i,j-1))+\Psi \Delta t .
\end{array}
\end{equation}

Next, we denote the error introduced by approximating $n(k,i,j)$ by $u(k,i,j)$  as
\begin{equation}
\label{(12)}
w(k,i,j)=n(k,i,j)-u(k,i,j)
\end{equation}
From \eqref{(6)} and \eqref{11}, we obtain
\begin{equation}\label{13}
\begin{array}{ll}
w(k+1,i,j)=&\left(1-v\frac{\Delta t}{\Delta a}-2D\frac{%
\Delta t}{(\Delta x)^{2}}\right) w(k,i,j)+v\frac{\Delta t}{\Delta a} w(k,i-1,j)
\vspace{2mm}\\
&+D\frac{\Delta t}{(\Delta x)^{2}}(w(k,i,j+1)+w(k,i,j-1))+\Psi \Delta t.
\end{array}
\end{equation}

From the stability condition \eqref{14} and using \eqref{13}, we obtain
\begin{equation}\label{15}
\begin{array}{ll}
|w(k+1,i,j)| \leq &\left( 1-v\frac{\Delta t}{\Delta a}-2D\frac{\Delta t}{(\Delta x)^{2}}\right)
|w(k,i,j)| +v\frac{\Delta t}{\Delta a}|w(k,i-1,j)| + \vspace{2mm}\\
& +D\frac{\Delta t}{(\Delta x)^{2}}(|w(k,i,j+1)|+|w(k,i,j-1)|)+|\Psi|
\Delta t
\end{array}
\end{equation}

For every \textit{k} we define
$$
M_{k}=\max\limits_{i,j} |w(i,j,k)|.
$$

From \eqref{15}, we obtain successively 
\[
M_{k+1}\leq \left( 1-v\frac{\Delta t}{\Delta a}-2D\frac{%
\Delta t}{(\Delta x)^{2}}\right) M_{k}+v\frac{\Delta t}{\Delta a}M_{k}+2D%
\frac{\Delta t}{(\Delta x)^{2}}M_{k}+\left| \Psi \right| \Delta t ,
\]

\[
M_{k+1}\leq \ M_{k}+|\Psi| \Delta t.
\]
Summing over $k=0,1,...,m-1,\,\,m\leq K$

\[
\sum_{k=0}^{m-1}\
M_{k+1}\leq \sum_{k=0}^{m-1}\ M_{k}+m\Delta t |\Psi| ,
\]
\begin{equation}\label{16}
M_{m}\leq  M_{0}+m\Delta t |\Psi| .
\end{equation}

Let us observe that $u(0,i\Delta a,j\Delta x)=n(0,i\Delta a,j\Delta x)$, and therefore we have

\[
w(0,i,j)=w(0,i\Delta a,j\Delta x)=0, \ \ \ (i=0,1,...,I\ ,\ j=0,1,...,J).
\]

Bearing in mind that $M_{0}=0$ and $\Delta t=\frac{T}{K}$, from (\ref{16}) we obtain

\[
M_{m}\leq m\frac{T}{K} |\Psi|.
\]

Since  $\frac{m}{K}\leq 1$, we also conclude that $M_{m}\leq  |\Psi| T$.\\

On the other hand, from \eqref{10}, we finally obtain
\[
M_{m}\leq \left( A\Delta t+B\Delta
a+C(\Delta x)^{2}\right) T,\ \ \ \ m=0,1,...,K.
\]

This proves that
\[
\max\limits_{k,i,j} |n(k,i,j)-u(k,i,j)|
\rightarrow 0 \quad \text{as} \quad \Delta t\rightarrow 0, \ \Delta
a\rightarrow 0, \ \Delta x\rightarrow 0.
\]
The previous reasoning is valid for both components of the solution of the system \eqref{(2)}, therefore
\[
\max\limits_{k,i,j} |n_{1}(k,i,j)-u_{1}(k,i,j)|
\rightarrow 0 \quad \text{and} \quad \max\limits_{k,i,j}
|n_{2}(k,i,j)-u_{2}(k,i,j)| \rightarrow 0,
\]
as  $\Delta t\rightarrow 0, \ \ \Delta a\rightarrow 0, \ \Delta x\rightarrow 0$.

\section{Dynamics of the simplified model from $E2$ to $E4$}

We now study the dynamics of the model between $E2$ and $E4$. It is assumed that only symmetric mitosis takes place during that period. Then, we use the model (discretized by finite differences) given by equations \eqref{(4)}.

\subsection{Experimental data recorded from the optic tectum at $E2$}

Records of number of cells counts performed on successive segments along the optic tectum cephalic-caudal axis were used as initial data for the model. Such information at day $E2$ corresponds to \textit{16 segments, 25}$\mu m$. We have named these segments ``sections''. Every segment is identified by a value of the abscissa $x$.

The initial model length is $L_{E2}=16\times 25\mu m=400\mu m$ ,
and the total initial cells number is $N_{E2}=800$.

The distribution of cells number per section is reported in TABLE 1 (Annexe)

\subsection{Determination of parameters}

{\bf Cell cycle duration between $E2$ and $E4$.} According to experimental data, in the interval $E2-E4$, the number of cells increases from approximately around 800 to around 9200 cells. Assuming that during this period the whole neuro-epithelial cells population exhibits a homogeneous proliferative behavior with a cell cycle duration around a mean value, then the averaged duration can be estimated as 13.6 hours.
\\
$\bullet$ Cell cycle duration: $G_{1}+S+G_{2}+M = 13h36min = 13.6 h.$
\\
Intermitotic cycle duration $\ G_{1}+S+G_{2}$ \ is 96 \% the cell cycle duration; hence,
\\
$\bullet$ Duration of $G_{1}+S+G_{2} = 13.056 h \cong  13h 03min$.

The mitotic cycle duration $M$ represents $4\%$ the cell cycle duration :
\\
$\bullet$ Mitotic duration $M = 0.544 h = 33 min$.
\\

\noindent {\bf Number of cycles $G_{1}+S+G_{2}+M$ between E2 and E4.}

\[
n_{E_{2-4}}=\frac{48hs}{13.6hs}\cong 3.53
\]

\noindent {\bf  Total number of cells at E4.} Using the formula
\[
N_{E4}=N_{E2}\times 2^{n_{E_{2-4}}}
\]
we obtain
\[
N_{E4}=800\times 2^{3.53}\cong 9241
\]

\noindent {\bf Cellular density.} It was observed experimentally that the average cellular density $\delta
_{E_{2-4}}$ between $E2$ and $E4$ remained approximately constant,
$3.85\frac{cel}{\mu m}$. Then, it will be assumed that between $E2$ and $E4$ the model average cell density corresponds to this value.

\[
\delta _{E_{2-4}} = 3.85\frac{cel}{\mu m}
\]

\noindent {\bf  Model length at $E4$.}
\[
L_{E4}=\frac{%
N_{E4}}{\delta _{E_{2-4}}}=\frac{9241cel}{3.85\frac{cel}{\mu m}}\cong
2400\mu m
\]

\noindent {\bf  Number of $25\mu m$  thickness sections comprising the
model at $E4$}.
\[
\ D_{E_{2-4}}\cong \frac{2400\mu m}{400\mu m}=6
\]

\noindent {\bf  Calculation of $v_{1}$ and $v_{2}$ between $E2$ and
$E4$.} 
The following formulae define $v_{1}$and $v_{2}$,
\\
$\bullet$ 
$\frac{1}{v_{1}}=n^{er}$ of interphases $G_{1}+S+G_{2}$ per day $=\frac{24hs}{interphase\text{ }duration\text{ }in
\text{ }hs}$
\\
$\bullet$ 
$\frac{1}{v_{2}}=n^{er}$ of mitosis $M$ by day $=
\frac{24hs}{mitotic\text{ }duration\text{ }in\text{ }hs}$
\\

In our case we have

\[
\frac{1}{v_{1}}\cong \frac{24}{13}\cong 1.85 \ \ \ \ \frac{1}{v_{2}}\cong \frac{24}{0.544}\cong 44.12
\]
and we obtain
\[
v_{1}\cong 0.54 , \qquad \qquad  v_{2}\cong 0.023 .
\]

\noindent {\bf Time interval.} 
We take $t=0$ at $E2$ and \ $t=T=48hs$ at $E4$.
We decompose the interval in sub-intervals of $\ 10$ \ minutes and we take
$\Delta t=10min=\frac{1}{6}hour$.\\
The index will be \ $k=0,1,...,48\times 6=288$ \ , $K=288.$
\\

\noindent {\bf  Interval on the longitudinal axis.}
We decompose the interval $[0,400]$ in 16 sub-intervals of length $\Delta x=25\mu m$.\\
The index will be \ $j=0,1,...,16$ \ \ \ $J=16$.\\

\noindent {\bf  Interval of variation of $a$.} 
The interval of variation of $a$ is $[0,1]$.
We will choose $\Delta a$ \ to satisfy the stability condition
\[
v_{1}\frac{\Delta t}{\Delta a}+2D\frac{\Delta t}{%
(\Delta x)^{2}}\leq 1,\ \ \ \ \ v_{2}\frac{\Delta t}{\Delta a}+2D%
\frac{\Delta t}{(\Delta x)^{2}}\leq 1.
\]

Since
\begin{equation}\label{17}
v_{1}=0.54\frac{1}{h}, \ \Delta t=\frac{1}{6}h, \  D=6
\frac{(\mu m)^{2}}{h}, \ \Delta x=25\mu m,  \ v_{2}=0.023\frac{1}{h}
\end{equation}
the condition of stability will be the following one:
\[
0.54\frac{1/6}{\Delta a}+12\frac{1/6}{625}\leq 1, \ \
\ \ \ \ 0.023\frac{1/6}{\Delta a}+12\frac{1/6}{625}\leq 1.
\]
From the stability condition \eqref{14}, it is sufficient to choose $\Delta a$ such that
\[
0.54\frac{1/6}{\Delta a}+12\frac{1/6}{625}\leq 1,\ \ \
or\ \ \ \Delta a>0.09.
\]
We will adopt that index will be  $i=0,1,...,10$ that is 
\begin{equation}\label{18}
I=10, \qquad \qquad \Delta a=0.1\ . 
\end{equation}

\subsection{Induction formula for  $u_{1}$}

\begin{equation}\label{19}
\begin{array}{ll}
u_{1}(k+1,i,j)& = u_{1}(k,i,j)-v_{1}\frac{\Delta t}{\Delta a}
(u_{1}(k,i,j)-u_{1}(k,i-1,j))+
\vspace{2mm}\\
&\quad +D\frac{\Delta t}{(\Delta x)^{2}}(u_{1}(k,i,j+1)-2\,u_{1}(k,i,j)+u_{1}(k,i,j-1)).
\end{array}
\end{equation}
Summing over $i$ from $i=0$ \ to $i=I$ \ in (19), we obtain\ \
\begin{equation}\label{20}
\begin{array}{l}
 \sum\limits_{i=0}^{I}u_{1}(k+1,i,j)=\sum\limits_{i=0}^{I}u_{1}(k,i,j)-v_{1}\frac{\Delta
t}{\Delta a}\left\{
\sum\limits_{i=0}^{I}u_{1}(k,i,j)-\sum\limits_{i=0}^{I}u_{1}(k,i-1,j)\right\} +\vspace{2mm}\\
\hspace{5mm}+D\frac{\Delta t}{(\Delta x)^{2}}\left\{
\sum\limits_{i=0}^{I}u_{1}(k,i,j+1)-2\sum\limits_{i=0}^{I}\,u_{1}(k,i,j)+%
\sum\limits_{i=0}^{I}u_{1}(k,i,j-1)\right\}.
\end{array}
\end{equation}
Observe that
\[
\sum_{i=0}^{I}u_{1}(k,i,j)-%
\sum_{i=0}^{I}u_{1}(k,i-1,j)=u_{1}(k,I,j)-u_{1}(k,-1,j),
\]
therefore, from (\ref{20}) we obtain
\begin{equation}\label{21}
\begin{array}{l}
\sum\limits_{i=0}^{I}u_{1}(k+1,i,j)=\sum\limits_{i=0}^{I}u_{1}(k,i,j)-v_{1}\frac{\Delta t}{%
\Delta a}\left\{ u_{1}(k,I,j)-u_{1}(k,-1,j)\right\} + \\
\hspace{5mm} +D\frac{\Delta t}{(\Delta x)^{2}}\left\{
\sum\limits_{i=0}^{I}u_{1}(k,i,j+1)-2\sum\limits_{i=0}^{I}\,u_{1}(k,i,j)+%
\sum\limits_{i=0}^{I}u_{1}(k,i,j-1)\right\}.
\end{array}
\end{equation}
Let us also observe that the sum
\[U_{1}(l,r)=\sum_{i=0}^{I}u_{1}(l,i,r)
\]
is the total number of interphasic cells $(G1+S+G2)$ at instant $l$ and at section $r$.\\
Replacing in (\ref{21}) we obtain
\begin{equation}\label{22}
\begin{array}{ll}
U_{1}(k+1,j)&=U_{1}(k,j)-v_{1}\frac{\Delta t}{\Delta a}%
\left\{ u_{1}(k,I,j)-u_{1}(k,-1,j)\right\} +\vspace{2mm}\\
& +D\frac{\Delta t}{(\Delta x)^{2}}\left\{
U_{1}(k,j+1)-2\,U_{1}(k,j)+U_{1}(k,j-1)\right\}.
\end{array}
\end{equation}
$u_{1}(k,I,j)$ is the number of cells that complete the interphase $G1+S+G2$
at instant $k$, section $j$ .\\
$u_{1}(k,0,j)$ is the number of cells that begin the interphase $G_{1}+S+G_{2}$ at instant $k$, section $j$.

The number $u_{1}(k,0,j)-u_{1}(k,0,j)$ depends on the total number of
intephasic cells $G_{1}+S+G_{2}$, at instant $k$ and section $j$ ; therefore, it is a function of $U_{1}(k,j)$.\\
It is difficult establishing this dependence law. It will be assumed that there exist a constant $\alpha (k,j)$ , that depends on instant
$k$ and section $j$, such that
\begin{equation}\label{23}
u_{1}(k,0,j)-u_{1}(k,0,j)=\alpha (k,j)\ U_{1}(k,j).
\end{equation}
Now we replace in (\ref{22}), bearing in mind (\ref{23}) and the values of $\
v_{1},\Delta t$ , $\Delta x$ , $\Delta a$ , $D$ in (\ref{17}), to obtain
\begin{equation}\label{24}
U_{1}(k+1,j)=(1+0.9\alpha (j))
U_{1}(k,j)+0.0032\left( \tfrac{U_{1}(k,j+1)+U_{1}(k,j-1)}{2}
-U_{1}(k,j)\right).
\end{equation}

\subsection{Approximation of the induction formul\ae}

Since, for each $k$, the function $U_{1}(k,j)$ is continuous in $j$, we
conclude that
\[
0.0032\left( \frac{U_{1}(k,j+1)+U_{1}(k,j-1)}{2}-U_{1}(k,j)\right)
\]
will be very small. Therefore, it will be ignored in a first approximation.
Then, the formula (\ref{24}) takes the following approximate form
\begin{equation}\label{25}
U_{1}(k+1,j)=\left(1+0.9\alpha (j)\right) U_{1}(k,j).
\end{equation}
Multiplying expressions (\ref{25}) for $k=0,1,...,287$ we obtain
\begin{equation}\label{26}
\begin{array}{rcl}
\prod\limits_{k=0}^{287}U_{1}(k+1,j)&=&(1+0.9\alpha (j))^{288}\prod\limits_{k=0}^{287}U_{1}(k,j),
\ \ \ or \vspace{2mm}\\
U_{1}(48hs,j)&=&(1+0.9\alpha (j))^{288}U_{1}(0,j).
\end{array}
\end{equation}
Now,$\ U_{1}(48hs,j)$ represents the number of cells generated between $E2$ and $E4$ by the cells that at $E2$ were located at segment $j$. Given that during this period the model length increases 6 fold, it is assumed that, at $E4$ these newly generated cells will expand over six $25\mu m$\ thickness segments.\\

Thus, cells located at segment at $E2$ give origin to cells that at $E4$ occupy
$s$ - sections from $s=1+6(j-1)$\ to $s=6j$ . The values of $j$\ will be
$j=1,2,...,16$, corresponding to 16 sections of the model at $E2$.

$U_{1}(48hs,s)$ indicates the number of cells located at segment of abscissa $s$ at $E4$.

$\bar{U}_{1}(48hs,s)_{s=1+6(j-1)}^{s=6j}$ \ will indicate the average number of cells in the segments of abscissa $s$ at $E4$, for $1+6(j-1)\leq s<6j$.\\
We have
\[
\bar{U}_{1}(48hs,s)_{s=1+6(j-1)}^{s=6j}=\allowbreak \frac{U_{1}(48hs,j)}{6}.
\]

For the sake of brevity, we will sometimes write $\bar{U}_{1}(48hs,s)$ \ instead of $\bar{U}_{1}(48hs,s)_{s=1+6(j-1)}^{s=6j}$ for the corresponding value of $j$ which appears in the above mentioned formula. We will have
\begin{equation}\label{27}
\begin{array}{ll}
\bar{U}_{1}(48hs,s)=\frac{1}{6}(1+0.9\alpha(j))^{288}U_{1}(0,j),& 1+6(j-1)\leq s<6j,\ or \vspace{2mm}\\
\bar{U}_{1}(48hs,s)=\beta (j)U_{1}(0,j),&j=1,2,...,16,
\end{array}
\end{equation}
where
\[
\beta (j)=\frac{1}{6}\left( 1+0.9\alpha (j)\right) ^{288}.
\]
We define the cell proliferation at segment $j$\ as the quotient
\[
cel.prol.(j)=\frac{U_{1}(48hs,j)}{U_{1}(0hs,j)}.
\]
Since $U_{1}(48hs,j)=6$ $\bar{U}_{1}(48hs,s)$, we obtain
\[
cel.prol.(j)=6 \beta(j).
\]
$\beta (j)$  is a coefficient proportional to the cell proliferation at segment $j$ .\\
The relationship between $\alpha (j)$ and $\beta (j)$ for $j=0,1,...,16$ is
\[
\alpha (j)=\frac{1}{0.9}\left( [6\beta (j)]^{\frac{1}{288}}-1\right)
 \ \ \ \ \beta (j)=\frac{1}{6}\left( 1+0.9\alpha (j)\right) ^{288}
\]
It will be assumed that proliferation is constant along the model, with a value equal to an average proliferation. \ This is equivalent to assume that
$\ \ \alpha (j)=\alpha $\ \ \ and $\ \beta (j)=\beta $ \ \ are independent
of $j$.\\

Finally, we will have
\begin{equation}\label{28}
\bar{U}_{1}(48hs,s)=\beta(j) U_{1}(0,j)\ \ \ \ \ \ \ \ j=1,2,...,16.
\end{equation}
We will calculate the constants $\alpha $ and $\beta $ . We will consider a hypothetical segment $j$ in which the number of cells will be obtained by calculating the average of the total number of cells of 16 sections in which the model was divided at $E2$. This number is
\[
\frac{800}{16}=50.
\]
Likewise, we will consider that cells of section $j$ originate a number of cells that at $E4$ occupy 6 sections of the model. The number of cells at each one of the above mentioned sections is obtained by calculating the average of the totality of cells of the 96 sections comprising the model. This average value is
\[
\frac{9241}{96}\cong 96.26
\]
and, therefore, 578 cells correspond to 6 sections.
We will apply the formula (\ref{26}) taking
\[
U_{1}(0,j)=50 \ \ \ \ \ \ \ \ \ \ \ \ \ \ U_{1}(48hs,j)=578.
\]
We obtain \ \ \ \ \ \ \ \ \ \ \ \ \ \ \ \ $578=(1+0.9\alpha)^{288} 50$\,\\
and therefore
\[
\alpha =\frac{10}{9}\left( \left[\frac{578}{50}\right]^{\frac{1}{288}}-1\right).
\]
We conclude that
\begin{equation}\label{29}
\alpha \cong 0.0095 \ \ \ \ and \ \ \ \ \beta \cong 1.93
\end{equation}
Replacing in (\ref{28}) we obtain
\begin{equation}\label{30}
 \bar{U}_{1}(48hs,s)=1.93 U_{1}(0,j)
\ \ \ 1+6(j-1)\leq s\leq 6j , \ j=1,2,...,16.
\end{equation}

\subsection{Number of cells given by formula (30) and comparison with experimental data}

The distribution of cells number per section is reported in TABLE 2 (Annexe)

\subsection{Approximation for linear interpolation of the number of cells in the intermediate sections and comparison with experimental data corresponding to $E4$}

We will put
\[
\bar{U}_{1}(s)= \bar{U}_{1}(48hs,s), \ \ \ \ 1+6(j-1)\leq s\leq 6j , \ j=1,2,...,16.
\]

We will use the formula of linear interpolation
\begin{equation}\label{31}
\bar{U}_{1}(s)=\bar{U}_{1}(1+6(j-1))+\frac{\bar{U}
_{1}(1+6j)-\bar{U}_{1}(1+6(j-1))}{6}(s-1)
\end{equation}
for $1+6(j-1)\leq s\leq 6j$ ,$j=1,2,...,16$.\\
The value $\bar{U}_{1}(97)$\ used in the calculation for the slope in
interval $91\leq s\leq 96$ is obtained by means of the following trick. Let us observe that in the initial tip of the model at $E2$ the variation in the cell number between segments $j=1$ and $j=2$ is, in increasing sense
\[
U_{1}(0,j=2)-U_{1}(0,j=1) =7
\]
This induces us to extrapolate TABLE 2 by putting
$U_{1}(0,j=17)=42$,
after which the variation at the model final tip will be in diminishing
sense, of the same magnitude,
\[
U_{1}(0,j=17) -U_{1}(0,j=16)=-7
\]
Applying formula (\ref{30}) for $j=17$ we obtain
\[
\bar{U}_{1}(s)=\bar{U}_{1}(48hs,s)\cong
1.93\times 42\cong 81 \ \ \ for \ \ 97\leq s\leq 102.
\]
We take the value $\bar{U}_{1}(97)=81$.

\begin{figure}[ht]
\begin{center}
\includegraphics[width=12cm, height=6cm]{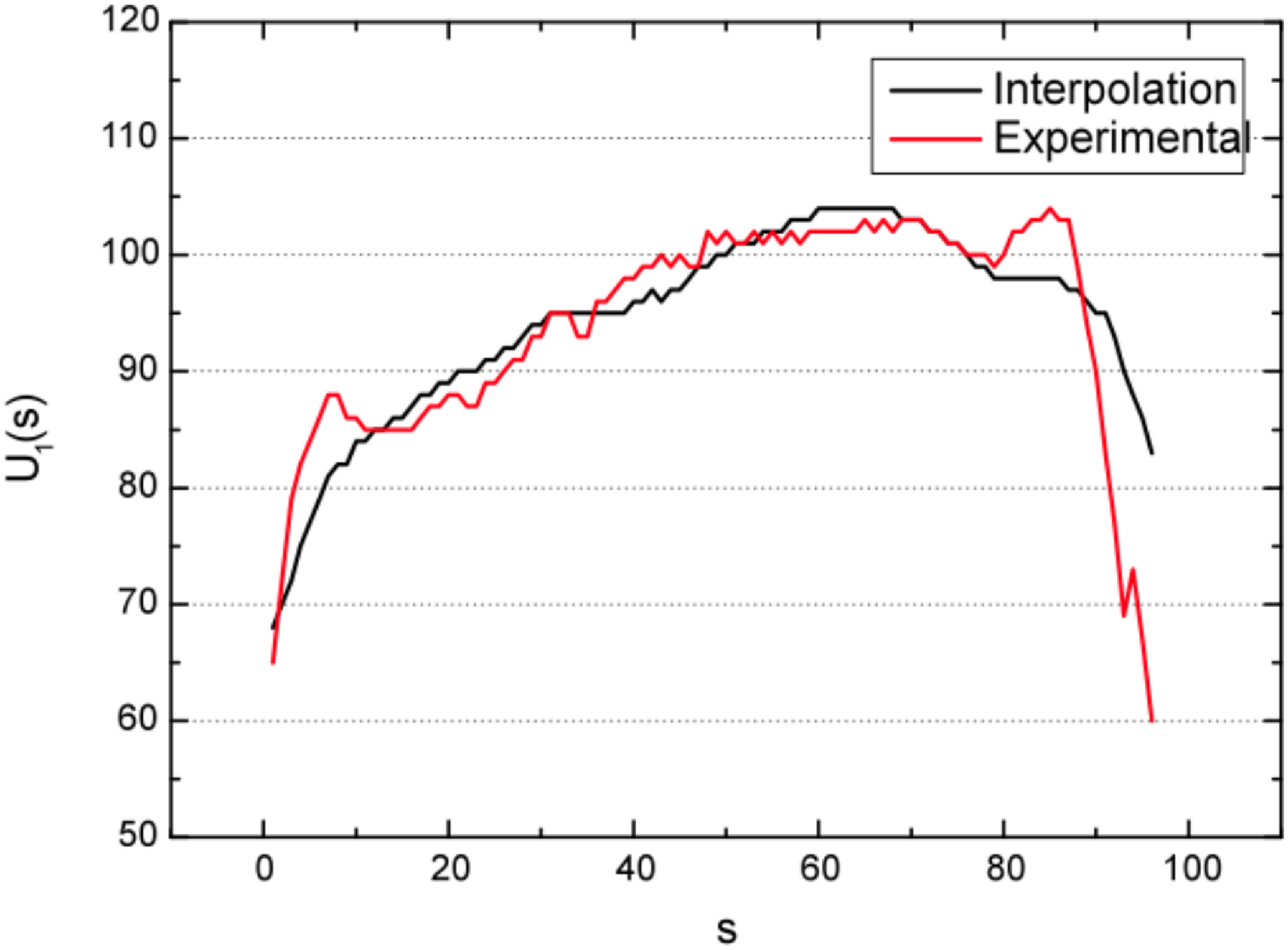} 
 \end{center}
\vspace{-.6cm} 
\caption{\label{figure1} 
Total number of cells per section $s$ in E4. 
} \end{figure}

Total number of cells per section $s$ at $E4$ is reported in TABLE\ 3 (Annexe) and Fig. \ref{figure1}.

\section{Dynamics of the model between $E4$ and $E6$}

\subsection{Initial data}

The number of cells calculated by the model for an $E4$ optic tectum was taken as initial data to calculate the evolution between $E4$ and $E6$ and to estimate the number of cell each section will originate over that period.\\

At $E$4 the model is composed of 96 segments , $25\mu m$ thickness.\\
The model initial length at $E4$ is $L_{E4}=96\times 25\mu m=2400\mu m$\\
The total initial number of cells of the model at $E4$ is $N_{E4}=9241$.

Table 3 illustrates spatially ordered distribution of number of cells/segment mathematically (left) and empirically (right) estimated.

\subsection{Determination of parameters}

\noindent {\bf Cell cycle duration between $E4$ and $E6$.} 
Based on assumptions introduced in a previous section (4.2.) the average duration of the complete cell cycle $G_{1}+S+G_{2}+M$ can be estimated as around as around 38 hs 15 min.\\
Cell cycle duration $G_{1}+S+G_{2}+M = 38.25 h =38 h 15 min$.\\
Assuming that the interphase $G_{1}+S+G_{2}$ is 96 \% the complete cycle,\\
Interphase duration $G_{1}+S+G_{2} = 36.72 h \cong  36 h 43 min$.\\

The mitotic length $M$ is 4 \% the complete cycle length,\\
Mitotic duration $M = 4 \% \  {\rm of} \ 38.25 h = 1.53 h \cong 1h 32 min\ $.\\

\noindent {\bf Number of cycles $G_{1}+S+G_{2}+M$ between $E4$ and $E6$.} 
\[
n_{E_{4-6}}=\frac{48hs}{38.25hs}\cong 1.255
\]

\noindent {\bf  Total number of cells of the model at $E6$.} 
Using the formula
\[
N_{E6}=N_{E4}\times 2^{n_{E_{4-6}}}
\]
we obtain
\[
N_{E6}=9241\times 2^{1.255}\cong 22055
\]

\noindent {\bf  Cellular density.} It was experimentally observed that the average cell density $\delta
_{E_{4-6}}$ of the Tectum between $E4$ and $E6$ remains approximately constant,
with a value of $3.85\ cel/\mu m$.

We will suppose therefore that $\delta _{E_{4-6}}=3.85\ cel/\mu m$.\\

\noindent {\bf  Model length at $E6$}. 
We have
$L_{E6}=\frac{N_{E6}}{\delta _{E_{4-6}}}=\frac{22055\,cel\text{\ }}{3.85cel/\mu m}\cong 5729\ \mu m$\\

\noindent {\bf  Number of $25\mu m$ thickness segments comprising the model at $E6$.} 
\\
$n^{er}$ of segments =$\frac{L_{E6}}{25\mu m}=\frac{5729\mu m}{25\mu m}\cong 230$\\

\noindent {\bf  Coefficient of cell diffusion between $E4$ and $E6$.} \\
$D_{E_{4-6}}=\frac{5729\mu m}{2400\mu m}=2.387...\cong \ 2.4 \mu m$\\

\noindent {\bf  Calculation of $v_{1}$ and $v_{2}$ between $E4$ and $E6$.} 
The formulae that define $v_{1}$ and $v_{2}$ are the following ones
\\
$\bullet$ $\frac{1}{v_{1}}=n^{er}$ of interphases $G_{1}+S+G_{2}$ per day  
$=\frac{24hs}{interphase\text{ }duration\text{ }in\text{ }hs}$
\\
$\bullet$ $\frac{1}{v_{2}}=n^{er}$ of mitosis $M$ per day = $\frac{24hs}{%
mitosis\text{ }duration\text{ }in\text{ }hs}$\\

We have  $\frac{1}{v_{1}}=\frac{24}{36.72}$ and $\frac{1}{v_{2}}=\frac{24}{1.53}$
therefore 
$$
v_{1}=\frac{36.72}{24}=1.53, \qquad \qquad v_{2}=\frac{1.53}{24}\cong 0.064.
$$

\noindent {\bf Time interval.} We take $t=0$\ at $E4$ and $t=T=48hs$  at $E6$. 
We decompose the interval in sub-intervals of 10 minutes each and we take
$$
\Delta t=10\min =\frac{1}{6}hour.
$$
The index will be \ \ $k=0,1,...,48\times 6=288,$  $K=288.$\\

\noindent {\bf  Interval on the longitudinal axis.} 
We decompose the interval $[0,2400]$ in 96 sub-intervals of length
$\Delta x=25\mu m$. The index will be \ \ $j=0,1,...,96$ \ , \ $J=96.$\\

\noindent {\bf  Interval of variation of $a$.}  The interval of variation for $a$ is still $[0,1]$. 
We choose $\Delta a$ so that the stability condition \eqref{14} is satisfied
\[
\nu _{1}\frac{\Delta t}{\Delta a}+2D\frac{\Delta t}{
(\Delta x)^{2}}\leq 1 \ \ \ \ \ \ \ \ \ \nu _{1}\frac{\Delta t}{\Delta a}
+2D\frac{\Delta t}{(\Delta x)^{2}}\leq 1.
\]
Since
\begin{equation}\label{32}
v_{1}=1.53\frac{1}{h},\, \Delta t=\frac{1}{6}h, \,
D=2.4\frac{(\mu m)^{2}}{h},\, \Delta x=25\mu m,\,  v_{2}=0.064\frac{1}{h}.
\end{equation}
the stability condition becomes 
\[
1.53\frac{1/6}{\Delta a}+4.8\frac{1/6}{625}
\leq 1, \qquad  0.064\frac{1/6}{\Delta a}+4.8\frac{1/6}{625}\leq 1,
\]
and it is sufficient to choose $\Delta a>0.256$. 
We adopt
\begin{equation}\label{33}
\Delta a=\frac{1}{3}
\end{equation}
The index will be \ $i=0,1,2,3$ \ \ \ \ \ $I=3$.

\subsection{Induction formula for $u_{1}$}

\begin{equation}\label{34}
\begin{array}{ll}
u_{1}(k+1,i,j)=&u_{1}(k,i,j)-v_{1}\frac{\Delta t%
}{\Delta a}(u_{1}(k,i,j)-u_{1}(k,i-1,j))+ \vspace{2mm}\\
& +D\frac{\Delta t}{(\Delta x)^{2}}(u_{1}(k,i,j+1)-2\,u_{1}(k,i,j)+u_{1}(k,i,j-1)).
\end{array}
\end{equation}

Proceeding along the same lines as in case of the model dynamics between $E2$ and $E4$, we obtain the formula
\begin{equation}\label{35}
\begin{array}{ll}
U_{1}(k+1,j)=&U_{1}(k,j)-v_{1}\frac{\Delta t}{\Delta a}\left\{
u_{1}(k,I,j)-u_{1}(k,-1,j)\right\} + \vspace{2mm}\\
&+D\frac{\Delta t}{(\Delta x)^{2}}\left\{
U_{1}(k,j+1)-2\,U_{1}(k,j)+U_{1}(k,j-1)\right\},
\end{array}
\end{equation}
where $U_{1}(l,r)=\sum_{i=0}^{I}u_{1}(l,i,r)$  is the total number of
interphasic cells $G_{1}+S+G_{2}$ at instant $l$ and at section $r$,\\
$u_{1}(k,I,j)$ is the number of cells at the end of the interphase $
G_{1}+S+G_{2}$ at instant $k$ , and at section $j$,
\\
$u_{1}(k,0,j)$ is the number of cells starting the
interphase $G_{1}+S+G_{2}$ at instant $k$ , and at section $j$.
\\

As we did it for the model between $E2$ and $E4$, we suppose that there exists a constant $\alpha (k,j)$ such that
\begin{equation}\label{36}
u_{1}(k,0,j)-u_{1}(k,0,j)=\alpha (k,j) U_{1}(k,j),
\end{equation}
and that between $E4$
and $E6$ $\alpha $ is also independent of $k$, i.e. $\alpha (k,j)=\alpha (j)$.\\

Replacing in (\ref{35}), and bearing in mind (\ref{36}), the values of $v_{1},\Delta t,\Delta
x,\Delta a,D$ in (\ref{32}) and (\ref{33}), we obtain
\begin{equation}\label{37}
U_{1}(k+1,j)=\left( 1+0.765\alpha (j)\right)
U_{1}(k,j)+0.00064\left( \tfrac{U_{1}(k,j+1)+U_{1}(k,j-1)}{2}%
-U_{1}(k,j)\right).
\end{equation}

\subsection{Approximation of the induction formula}

The term
\[
0.00064\left( \frac{U_{1}(k,j+1)+U_{1}(k,j-1)}{2}-U_{1}(k,j)\right)
\]
effectively acquires very small values which, in a first approximation, will
be ignored.\\

Formula (\ref{37}) takes the following approximate form
\begin{equation}\label{38}
U_{1}(k+1,j)=\left(1+0.765\alpha (j)\right) U_{1}(k,j).
\end{equation}
Multiplying expressions (\ref{37}) for $k=0,1,...,287$ we obtain
\begin{equation}\label{39}
U_{1}(48hs,j)=\left(1+0.765\alpha (j)\right) ^{288}U_{1}(0,j).
\end{equation}
$U_{1}(48hs,j)$ represents the number of cells generated between $E4$ and $E6$ by the cells that at $E4$ were located at section $j$ . Given that during this
period the model length undergoes a 2.4 fold increase, it is assumed that,
at $E6$ these newly generated cells will expand over 2.4 segments of $25\mu m$
thickness.\\

We will indicate with $s$ the abscissa of transverse sections of the model at $E6$. At each abscissa $s=1,2,...,230$ will correspond a segment of $25\mu
m $ thickness.\\

It is assumed that the cells located at the $j$-th segment at $E4$ give origin to cells that, at $E6$, occupy the part of the model corresponding to $2.4(j-1)\leq s<2.4$ $j$, $j=1,2,...,96$.\\

$U_{1}(48hs,s)$ indicates the number of cells located at the segment of
abscissa $s$ at $E6$.\\

$\bar{U}_{1}(48hs,s)_{s=2.4(j-1)}^{s=2.4j}$ will indicate the average number
of cells at the segment of abscissa $s$ at $E6$ , for $2.4(j-1)\leq s<2.4$ $j$.
We have
\begin{equation}\label{40}
\bar{U}_{1}(48hs,s)_{s=2.4(j-1)}^{s=2.4j}=\frac{U_{1}(48hs,j)}{2.4}
\end{equation}
The values of $j$ will be $j=1,2,...,96,$ which correspond to 96 segments of
the model at $E4$.\\

In order to abbreviate notations, in the following paragraphs, the formula
\[
\bar{U}_{1}(48hs,s)=\bar{U}_{1}(48hs,s)_{s=2.4(j-1)}^{s=2.4j}
\]
will simply be written as $\bar{U}_{1}(48hs,s)$ for the corresponding value of $j$.\\

We will have
\[
\bar{U}_{1}(48hs,s)=\frac{1}{2.4}\left( 1+0.765\alpha (j)\right) ^{288}U_{1}(0,j)
\]
for $2.4(j-1)\leq s<2.4$ $j$, $j=1,2,...,96$ \ or
\[
\bar{U}_{1}(48hs,s)=\beta (j) U_{1}(0,j)
\]
where
\begin{equation}\label{41}
\beta (j)=\frac{1}{2.4}\left( 1+0.765\alpha (j)\right) ^{288}
\end{equation}
We define the cell proliferation at segment $j$\ as the quotient
\[
prol(j)=\frac{U_{1}(48hs,j)}{U_{1}(0hs,j)}
\]
Since
\[
U_{1}(48hs,j)=2.4\ \bar{U}_{1}(48hs,s)
\]
we obtain
\[
prol(j)=2,4\beta (j).
\]
The $\beta (j)$ coefficient is proportional to the cell proliferation at
section $j$.\\

Experimental data indicate that during the interval $E4-E6$ cell proliferation is not uniformly distributed along the cephalic-caudal axis. Thus, a coefficient $%
\beta (j)$ is introduced by means of the following simplified expression
\begin{equation}\label{42}
\beta (j)=\left\{
\begin{array}{lcl}
0,1j+0.525&& if\  j=1,2,3,4,5,6 \\
1.125&& if\  j=6,7,...,48 \\
-0.01j+1.605&& if\ j=48,...,95,96
\end{array}
\right.
\end{equation}
From (\ref{41}) we obtain
\[
\alpha (j)=\frac{1}{0.765}\left( (2.4\beta (j))^{\frac{1}{288}%
}-1\right)
\]
and replacing by (\ref{42}) we obtain the final expression for $\alpha (j)$.\\

Finally, we will have
\begin{equation}\label{43}
\bar{U}_{1}(48hs,s)=
\beta (j) U_{1}(0,j) \ \ \ for \ \ \ j=1,2,...,96,
\end{equation}
where $\beta (j)$ is given by (\ref{42}).

\subsection{Number of cells given by formula (43) and comparison with experimental data}
Number of cells per section $s$ is reported in TABLE\ 4 (Annexe)\
\subsection{Approximation by linear interpolation of the cell number in
intermediate segments of the model at $E6$}
We will denote by $U_{1}(48hs,s)=U_{1}(s)$\ the number of cells at segment $s$ of the model at $E6$ obtained by linear interpolation between values of \ TABLE 4.\\

We will have
\begin{equation}\label{44}
U_{1}(s)=\left( 1-\frac{s-2.4(j-1)}{2.4}\right)
\overline{U}_{1}(2.4(j-1))+\frac{s-2.4(j-1)}{2.4}\overline{U}_{1}(2.4 j)
\end{equation}
for sections of the model in $E6$, with absciss\ae , $2.4(j-1)\leq
s<2.4\,j$ .\\

We apply formula (\ref{44}) to calculate $U_{1}(s)$\ for integer values of $s$.\\

\begin{figure}[ht]
\begin{center}
\includegraphics[width=12cm, height=6cm]{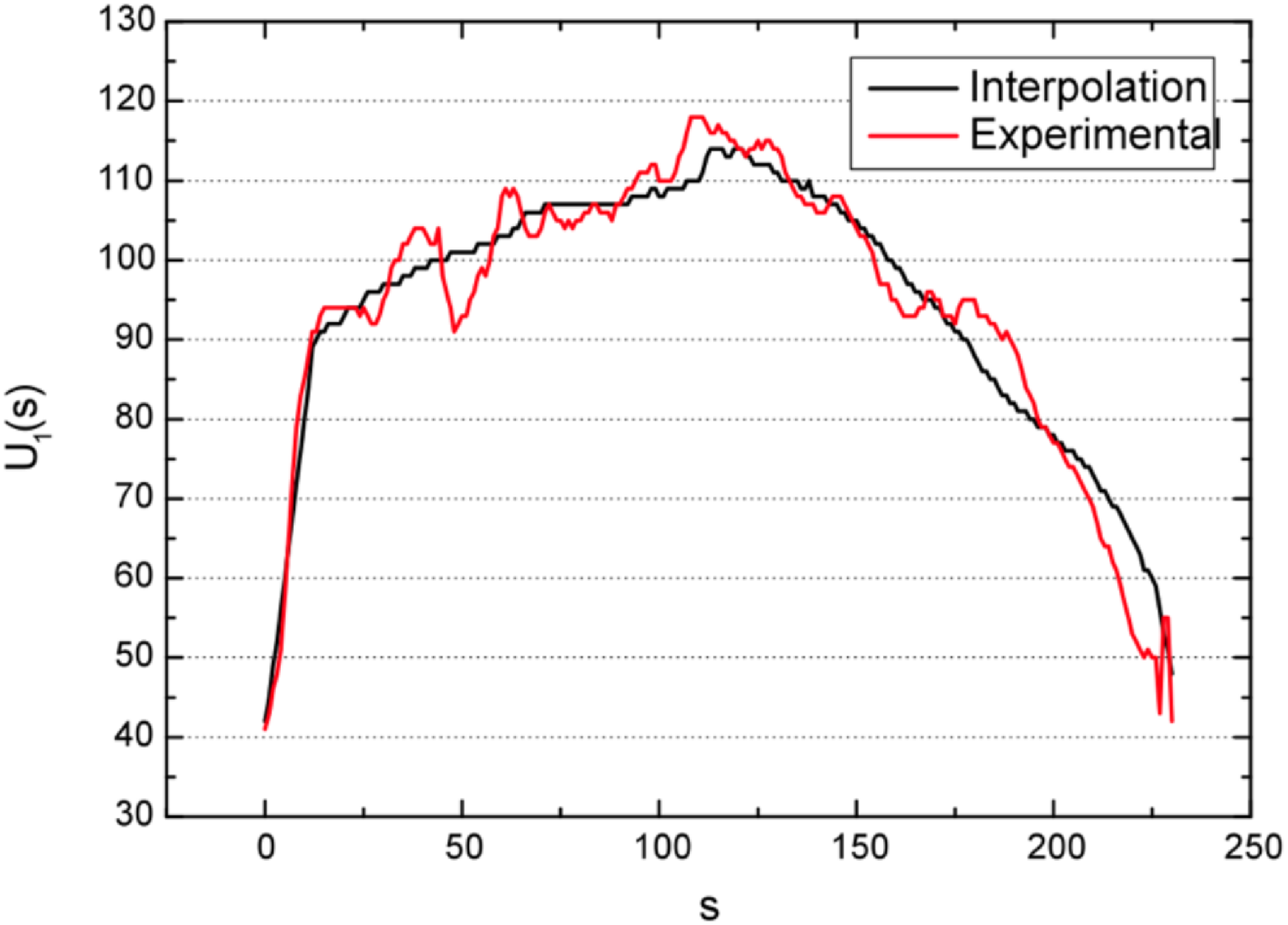} 
 \end{center}
\vspace{-.6cm} 
\caption{\label{figure2} 
Total number of cells per section $s$ in $E6$. 
} \end{figure}

Total number of cells per section $s$ in $E4$ is reported in TABLE\ 5
(Annexe)\ and Fig. \ref{figure2}.

\section{Discussion}

\subsection{How to analyze the developmental role of cell proliferation in a developing system?}

The proliferation of NE cells is not only involved in generating the appropriate number of cells in the CNS. A spatially organized cell proliferation activity also possesses morpho- and histogenetic effects \cite{[Cowan],[LaVail],[Alvarez],[Kahane]}.
These developmental effects should depend on the temporal dynamics and the spatial organization of the proliferative activity. Taking this notion in mind, the present paper attempts at designing a mathematical PDE-based model with the ability to describe the increase in cell number as a function of both time and space. The model emphasizes the occurrence of (a) changes in the cell cycle duration as function of time and (b) changes in the density of NE cells as a function of time and space.\\

In our opinion this sort of PDE-based models can be useful tools to clearly understand the cell proliferation spatial and temporal organization and to design formal models of cell proliferation that could account for its developmental effects.\\

This paper is a first step towards designing a model that could help to globally describe the temporal and spatial organization of the cell proliferation and to understand its specific roles in the genesis of CNS patterning and its supracellular complexity.\\

A current notion in theoretical biology is that the supracellular complexity is the net result of several DCBs operating simultaneously and interactively \cite{Gilbert}. The corticogenesis requires the integrated operation of several DCBs. Such organized operation requires spatially organized cell signaling processes which do not only control cell proliferation and migration but also establish cell-type differences and also determines tissue patterning \cite{[Freeman]}. It is currently accepted that these kinds of directive influences are mediated by gradients of specific signaling molecules, the so-called morphogens. A developmental gradient may be defined as an asymmetry in morphogen concentration with a maximum at the site of the morphogen-secreting cells and decaying as a function of the distance from that site \cite{[1]}.\\

As an example, the non-uniform distribution of a morphogen in an embryonic field differentially determines the fate and phenotype of those cells \cite{[Green]}. In a same way, it can be proposed that the non-uniform spatial distribution of developmentally active signals influencing the proliferative activity should result in a non-uniform distribution of the density of proliferating cells and also in the existence of space-dependent differences in cell density along a defined spatial axis.\\

\subsection{Stochastic approaches on the NE cell proliferation dynamics during the OT corticogenesis}

In the absence of directive influences a DCB should operate at random. In such condition a biological time or space signal representing the DCB operation or representing a quantifiable direct developmental effect should display a random white noise-like dynamics.\\

A previous work devoted to characterize the postmitotic neuronal migration dynamics in the developing OT by means of standardized methods of frequency signal analyses \cite{[Rapacioli01]}, proposes a mathematical expression that globally describes how cell migration operates in space. This expression includes deterministic and stochastic components. The deterministic component can be described as the sum of a linear global trend and a Fourier series while the stochastic component can be described as a non-correlated process (white noise).\\

It is reasonable considering that, in the absence of directive influences the NE cells proliferation should operate at random and, consequently, the intercalation of newly generated cells into the whole system should be uniformly distributed. Under such condition the OT growth should proceed as a uniformly distributed tangential expansion; a characteristic that does not occur during the normal OT development.
In a previous paper, aiming at determining whether the NE cell proliferation exhibits a spatial organization, spatial records of proliferating (mitotic) NE cells were analyzed within the framework provided by the stochastic point process model \cite{[Mazzeo04a]}. These analyses showed that the NE cell proliferation records can also be described as comprising two components: (a) a non-stationary one that can be deterministically described as a trend plus (b) a stationary non-correlated stochastic point process. The deterministic trend can be interpreted as manifestation of an asymmetrically distributed controlling influence that installs space-dependent differences in the NE cells proliferation rate along the cephalic-caudal axis.\\

A more recent paper \cite{[Mazzeo04b]} indicates that the dynamics of a detrended NE cell proliferation record does not accurately coincide with a non-correlated stochastic point process with uniform probability distribution of intermitotic intervals. By contrast, stochastic point processes simulated with intermitotic intervals displaying lognormal, exponential and gamma distribution better approximate the NE cell proliferation records than those simulated with a uniform probability distribution.\\

\subsection{The present conclusions on the temporal and spatial organization of the NE cells proliferation during the OT corticogenesis. Morphogenetic correlations}

The PDE-based model presented in this paper was specifically design to obtain a better and more realistic characterization of the deterministic behavior exhibited by the proliferating NE cells. In our opinion, the present results appropriately describe the spatial organization of the proliferation activity and explain the differential growth the developing OT displays between $E2$Ð$E6$. The model accurately estimates the increase in the total number of cells along the defined temporal windows ($E2$-$E4$ and $E4$-$E6$). These increases closely coincide with values of cell density empirically estimated by means of cells counts performed on histological sections of developing OTs.\\

With regards to the spatial organization of the NE cells proliferation, the model appropriately reproduces the empirically recorded position-dependent differences in cell density observed in successively ordered $25\mu m$ length windows along the OT chepalic-caudal axis. In fact, figures 1 and 2 show that the cell density (number of cells / each $25\mu m$ windows) significantly changes as a function of their position along the cephalic-caudal axis. A comparison between these figures shows that, considering the OT as a whole, the NE cells proliferation significantly decreases as a function of the time (compare the profile of cell density values obtained at the end of the interval $E2$-$E4$ with the profile of cell density at the end of the interval $E4$-$E6$).\\

It is interesting remarking that, apart from reproducing the increase in the total number of cells, the profiles of cell density as a function of the space also allows accounting for some morphogenetic events. In fact, between $E2$Ð$E6$ the OT undergoes a remarkable change in shape and size. The differences in cell proliferation as a function of their position along the cephalic-caudal axis should result in a differential intercalation of new cells within the OT wall along that spatial axis. The addition of fewer new cells at both extremes of the OT longitudinal axis should imply a lower OT growth rate at those sites. It must be noted that the low rate of NE cell proliferation at both ends of the cephalic-caudal axis reveals a low NE proliferation rate along the dorsal midline. This lower proliferation activity should also produce a lower growth rate along the entire dorsal midline. This assumption is corroborated by the fact that between $E2$Ð$E6$ a deep medial groove appears along the dorsal midline. This change is a relevant morphogenetic effect since the appearance of a medial groove leads to the constitution of the right and the left OT hemispheres.\\

All these data about (a) the changes in the number of cells as a function of the time, (b) the space-dependent differences in NE cell proliferation along the cephalic-caudal axis together with (c) the adequacy between the changes in size and shape that can be inferred from the model results, strongly suggest that the PDE-based model here presented can be appropriately and advantageously used to more precisely describe and to more deeply understand the developmental effects of a spatially and temporally organized NE cell proliferation during the OT corticogenesis.



\newpage
\begin{center}
{\LARGE  ANNEXE}

\bigskip
\bigskip
\bigskip

{\Large TABLE \ 1}

\bigskip

\begin{tabular}{ccc}
Section & Number intermitotic cells & Number of mitotic cells \\
$j=$1 & 35 & 2 \\
2 & 42 & 1 \\
3 & 44 & 1 \\
4 & 46 & 1 \\
5 & 47 & 1 \\
6 & 49 & 1 \\
7 & 49 & 1 \\
8 & 50 & 2 \\
9 & 52 & 1 \\
10 & 53 & 1 \\
11 & 54 & 1 \\
12 & 54 & 1 \\
13 & 53 & 1 \\
14 & 51 & 1 \\
15 & 51 & 1 \\
16 & 49 & 1
\end{tabular}

\bigskip
\bigskip
\bigskip

{\Large TABLE\ 2}

\bigskip

\begin{tabular}{cccc}
$Section$ & $\underset{}{U_{1}(0,j)}$ & $\underset{according\text{ (30)}}{%
\bar{U}_{1}(48hs,s)}$ & $\underset{experimental\text{ }data}{\bar{U}%
_{1}(48hs,s)}$ \\
$j=1$ & 35 & 68 & 78 \\
2 & 42 & 81 & 86 \\
3 & 44 & 85 & 86 \\
4 & 46 & 89 & 87 \\
5 & 47 & 91 & 91 \\
6 & 49 & 95 & 94 \\
7 & 49 & 95 & 98 \\
8 & 50 & 96 & 100 \\
9 & 52 & 100 & 101 \\
10 & 53 & 102 & 101 \\
11 & 54 & 104 & 102 \\
12 & 54 & 104 & 103 \\
13 & 53 & 102 & 101 \\
14 & 51 & 98 & 101 \\
15 & 51 & 98 & 101 \\
16 & 49 & 95 & 72
\end{tabular}

\newpage

{\Large TABLE\ 3}

\bigskip

\begin{tabular}{ccc}
$\underset{order}{Section}$ & $\underset{\text{ }obtained\text{ }by\text{ }%
interpolation}{N^{er}cells\text{ }of\text{ }the\text{ }Model\text{ }at\text{
}E4}$ & $\underset{according\text{ }to\text{ }experimental\text{ }data}{%
N^{er}cells\text{ }of\text{ }the\text{ }Tectum\text{ }at\text{ }E4}$ \\
$s=1$ & 68 & 65 \\
2 & 70 & 72 \\
3 & 72 & 79 \\
4 & 75 & 82 \\
5 & 77 & 84 \\
6 & 79 & 86 \\
7 & 81  &88 \\
8 & 82 & 88 \\
9 & 82 & 86 \\
10 & 84 & 86 \\
11 & 84 & 85 \\
12 & 85 & 85 \\
13 &85 & 85  \\
14 & 86 & 85 \\
15 & 86 & 85 \\
16 & 87 & 85 \\
17 & 88 & 86 \\
18 & 88 & 87 \\
19 & 89 &87 \\
20 & 89 & 88 \\
21 & 90 & 88 \\
22 & 90 & 87 \\
23 & 90 & 87 \\
24 & 91 & 89 \\
25 & 91  & 89  \\
26 & 92 & 90 \\
27 & 92 & 91 \\
28 & 93 & 91 \\
29 & 94 & 93 \\
30 & 94 & 93 \\
31 & 95 &95  \\
32 & 95 & 95 \\
33 & 95 & 95 \\
34 & 95 & 93 \\
35 & 95 & 93 \\
36 & 95 & 96 \\
37 & 95 & 96 \\
38 & 95 & 97 \\
39 & 95 & 98 \\
40 & 96 & 98 \\
41 & 96 & 99 \\
42 & 97 & 99 \\
43& 96  & 100 \\
44 & 97 & 99 \\
45 & 97 & 100 \\
46 & 98 & 99 \\
47 & 99 & 99 \\
48 & 99 & 102\\
49 &100  &101 
\end{tabular}

\newpage

\begin{tabular}{ccc}
 $s=50$ & \ \qquad \qquad 100 \qquad \qquad \qquad & \ \qquad \qquad \qquad \qquad 102 \qquad \qquad \qquad \qquad \qquad \\
51 & 101 & 101 \\
52 & 101 & 101 \\
53 & 101 & 102 \\
54 & 102 & 101 \\
55 & 102 & 102 \\
56 & 102 & 101 \\
57 & 103 & 102 \\
58 & 103 & 101 \\
59 & 103 & 102 \\
60 & 104 & 102 \\
61 &104 &102 \\
62 & 104 & 102 \\
63 & 104 & 102 \\
64 & 104 & 102 \\
65 & 104 & 103 \\
66 & 104 & 102 \\
67 & 104 & 103  \\
68 & 104 & 102 \\
69 & 103 & 103 \\
70 & 103 & 103 \\
71 & 103 & 103 \\
72 & 102 & 102 \\
73 & 102 & 102 \\
74 & 101 & 101 \\
75 & 101 & 101 \\
76 & 100 & 100 \\
77 & 99 & 100 \\
78 & 99 & 100 \\
79 & 98  &99 \\
80 & 98 & 100 \\
81 & 98 & 102 \\
82 & 98 & 102 \\
83 & 98 & 103 \\
84 & 98 & 103 \\
85 & 98 & 104  \\
86 & 98 & 103 \\
87 & 97 & 103 \\
88 & 97 & 99 \\
89 & 96 & 94 \\
90 & 95 & 90 \\
91 & 95  &83 \\
92 & 93 & 77 \\
93 & 90 & 69 \\
94 & 88 & 73 \\
95 & 86 & 67 \\
96 & 83 & 60
\end{tabular}

\newpage 

{\Large TABLE\ 4}

\bigskip

\begin{tabular}{ccccc}
&  & $E4$ & $E6$ & $E6$ \\
&  & $U_{1}(0,j)$ & $\bar{U}_{1}(48hs,s)_{s=2.4(j-1)}^{s=2.4j}$ & $\bar{U}%
_{1}(48hs,s)_{s=2.4(j-1)}^{s=2.4j}$ \\
&  &  & $according$ $to$ $(43)$ & $experimental$ $data$ \\
$j=1$ & $\ 0\leq s<2.4$ \  & $68$ & $43$ & $43$ \\
$2$ & $\ 2.4\leq s<4.8$ \  & 70 & 51 & 48 \\
$3$ & $4.8\leq s<7.2$ & 72 & 59 & 62 \\
$4$ & $7.2\leq s<9.6$ & 75 & 69 & 78 \\
$5$ & $9.6\leq s<12$ & 77 & 79 & 86 \\
$6$ & $12\leq s<14.4$ & 79 & 89 & 91 \\
$7$ & $14.4\leq s<16.8$   & 81 &  91  &  94  \\
$8$ & $\ 2.4\leq s<4.8$ \  & 82 & 92 & 94 \\
$9$ & $4.8\leq s<7.2$ & 82 & 92 & 94 \\
$10$ & $7.2\leq s<9.6$ & 84 & 94 & 94 \\
$11$ & $9.6\leq s<12$ & 84 & 94 & 93 \\
$12$ & $12\leq s<14.4$ & 85 & 96 & 93 \\
$13$ & $28.8\leq s<31.2$ & 85 & 96  & 94 \\
$14$ & $\ 31.2\leq s<33.6$ \  & 86 & 97 & 98 \\
$15$ & $33.6\leq s<36$ & 86 & 97 & 101 \\
$16$ & $36\leq s<38.4$ & 87 & 98 & 103 \\
$17$ & $38.4\leq s<40.8$ & 88 & 99 & 104 \\
$18$ & $40.8\leq s<43.2$ & 88 & 99 & 102 \\
$19$ & $43.2\leq s<45.6$ &  89 & 100 & 100 \\
$20$ & $\ 45.6\leq s<48$ \  & 89 & 100 & 95 \\
$21$ & $48\leq s<50.4$ & 90 & 101 & 91 \\
$22$ & $50.4\leq s<52.8$ & 90 & 101   & 94 \\
$23$ & $52.8\leq s<55.2$ & 90 & 101  & 97 \\
$24$ & $55.2\leq s<57.6$ & 91 & 102  & 99 \\
$25$ & $57.6\leq s<60$ \  & 91  & 102 & 103 \\
$26$ & $ 60\leq s<62.4$ \  & 92 & 104   & 108 \\
$27$ & $62.4\leq s<64.8$ & 92 & 104  & 108 \\
$28$ & $64.8\leq s<67.2$ & 93 & 105  & 105 \\
$29$ & $67.2\leq s<69.6$ & 94 & 106  & 103 \\
$30$ & $69.6\leq s<72$ & 94 & 106  & 105 \\
$31$ & $72\leq s<74.4$ & 95& 107 &107\\
$32$ & $\ 74.4\leq s<76.8$ \  & 82 & 92 & 94 \\
$33$ & $76.8\leq s<79.2$ & 82 & 92 & 94 \\
$34$ & $79.2\leq s<81.6$ & 84 & 94 & 94 \\
$35$ & $81.6\leq s<84$ & 84 & 94 & 93 \\
$36$ & $84\leq s<86.4$ & 85 & 96 & 93 \\
$37$ & $86.4\leq s<88.8$ &  95& 107  & 106\\
$38$ & $\ 88.8\leq s<91.2$ \  & 95 & 107  & 107   \\
$39$ & $91.2\leq s<93.6$ & 95 & 107   & 109   \\
$40$ & $93.6\leq s<96$ & 96 & 108  & 110   \\
$41$ & $96\leq s<98.4$ & 96 & 108 & 111  \\
$42$ & $98.4\leq s<100.8$ & 97 & 109  & 111 \\
$43$ & $100.8\leq s<103.2$ & 96  &108  & 110 \\
$ 44$  & $\ 103.2\leq s<105.6$ \  & 97 & 109 & 111 \\
$45$ & $105.6\leq s<108$ & 97   & 109   & 115  \\
$46$ & $108\leq s<110.4$ & 98  & 110& 118  \\
$47$ & $110.4\leq s<112.8$ & 98 & 110  & 118  \\
$48$ & $112.8\leq s<115.2$ & 99 & 111  & 116 
\end{tabular}

\newpage
\begin{tabular}{ccccc}
\ \ \ $j=$49 & $115.2\leq s<117.6$ & \quad 100\quad  \ & \ \ \ \ \ \ $108$\ \ \
\ \ \ \ \ \ \ \ \ \  & $\ \ \ \ \ \ \ \ \ \ \ $110 \ \ \ \ \ \ \ \ \ \ \ \ \
\ \  \\
\ 50 & $\ 117.6\leq s<120$ \  &  100  & 109 \ \ \ \ \ \  & 111
\ \ \ \  \\
\ 51 & $120\leq s<122.4$  & 101  & 109 \ \ \ \ \ \  & 115 \ \
\ \  \\
\ 52 & $122.4\leq s<124.8$ & 101   & 110 \ \ \ \ \ \  & 118
\ \ \ \  \\
\ 53 & $124.8\leq s<127.2$ & 101  & 110 \ \ \ \ \ \  & 118 \ \
\ \  \\
\ 54 & $127.2\leq s<129.6$ &102  & 111 \ \ \ \ \ \  & 116 \
\ \ \ \\
\ 55 & $129.6\leq s<132$ & 102& $\ \ \ \ \ \ $108$\ \ $
\ \ \ \ \ \ \ \ \ \  & $\ \ \ \ \ \ \ \ \ \ \ $113 \ \ \ \ \ \ \ \ \ \ \ \ \
\ \  \\
\ 56 & $\ 132\leq s<134.4$ \  & 102 & 107 \ \ \ \ \ \  & 110
\ \ \ \  \\
\ 57 & $134.4\leq s<136.8$ & 103  & 107 \ \ \ \ \ \  & 108 \
\ \ \  \\
\ 58 & $136.8\leq s<139.2$ &103  & 106 \ \ \ \ \ \  & 107
\ \ \ \  \\
\ 59 & $139.2\leq s<141.6$ & 103 & 105 \ \ \ \ \ \  & 106 \ \
\ \  \\
\ 60 & $141.6\leq s<144$ & 104  & 105 \ \ \ \ \ \  & 107 \ \
\ \ \\
\ 61 & $144\leq s<146.4$ &104  & 102 \ \ \ \ \ \   & 108 \ \ \ \  \\
\ 62 & $\ 146.4\leq s<148.8$ \  & 104  & 102 \ \ \ \ \ \  & 107    \ \ \ \  \\
\ 63 & $148.8\leq s<151.2$ &104 & 101 \ \ \ \ \ \  & 105  \ \ \ \   \\
\ 64 & $151.2\leq s<153.6$ &104  & 100 \ \ \ \ \ \  & 103  \ \ \ \   \\
\ 65 & $153.6\leq s<156$ & 104 & 99 \ \ \ \ \ \  & 100 \ \ \ \   \\
\ 66 & $156\leq s<158.4$   & 104 & 98 \ \ \ \ \ \  & 97 \ \ \ \  \\
\ 67 & $158.4\leq s<160.8$   &104 & 97 \ \ \ \ \ \  & 95 \ \ \ \   \\
\ 68 & $\ 160.8\leq s<163.2$ &104 & 96 \ \ \ \ \ \  & 94 \ \ \ \   \\
\ 69 & $163.2\leq s<165.6$ & 103  & 94 \ \ \ \ \ \  & 93 \ \ \ \  \\
\ 70 & $165.6\leq s<168$ & 103     & 93 \ \ \ \ \ \  & 94 \ \  \ \  \\
\ 71 & $168\leq s<170.4$ & 103     & 92 \ \ \ \ \ \  &  96 \ \ \ \  \\
\ 72 & $170.4\leq s<172.8$ &102   & 90 \ \ \ \ \ \  & 95 \ \ \ \  \\
\ 73 & $172.8\leq s<175.2$     & 102  & 89 \ \ \ \ \ \ &  93 \ \ \ \  \\
\ 74 & $\ 175.2\leq s<177.6$\  &101   & 87 \ \ \ \ \ \  & 94 \ \  \ \ \\
\ 75 & $177.6\leq s<180$    &101 & 86 \ \ \ \ \ \  & 95 \ \  \ \ \\
\ 76 & $180\leq s<182.4$     &100 & 86 \ \ \ \ \ \  & 94 \ \ \ \ \\
\ 77 & $182.4\leq s<184.8$ & 99  & 83 \ \ \ \ \ \   & 92 \ \ \ \ \\
\ 78 & $184.8\leq s<187.2$ & 99 & 82 \ \ \ \ \ \    & 91 \ \ \ \ \\
\ 79 & $187.2\leq s<189.6$& 98 & 80  \ \ \ \ \ \  & 90 \ \ \ \ \\
\ 80 & $189.6\leq s<192$  & 98 & 79 \ \ \ \ \ \  & 89 \ \ \ \  \\
\ 81 & $192\leq s<194.4$   & 98  & 78 \ \ \ \ \ \  & 84 \ \ \ \  \\
\ 82 & $194.4\leq s<196.8$ & 98 & 77 \ \ \ \ \  \ & 81 \ \  \ \  \\
\ 83 & $196.8\leq s<199.2$ & 98 & 76 \ \ \ \ \ \ & 79 \ \  \ \  \\
\ 84 & $199.2\leq s<201.6$ & 98 & 75 \ \ \ \ \ \ & 78 \ \ \ \  \\
\ 85 & $201.6\leq s<204$ & 98  &   74 \ \ \ \ \ \  & 76  \ \ \ \ \\
\ 86 & $ 204\leq s<206.4$  & 98   & 73 \ \ \ \ \ & 74 \ \ \ \  \\
\ 87 & $206.4\leq s<208.8$ &  97  & 71 \ \ \ \ \ & 72  \ \ \ \  \\
\ 88 & $208.8\leq s<211.2$ &  97  & 70 \ \ \ \ \ & 69  \ \ \ \  \\
\ 89 & $211.2\leq s<213.6$ & 96 & 69 \ \ \ \ \ & 66 \  \ \ \  \\
\ 90 & $213.6\leq s<216$   & 95 & 67 \ \ \ \ \ & 63 \ \ \ \  \\
\ 91 & $216\leq s<218.4$    & 95 & 66 \ \ \ \ \ & 60  \ \ \ \  \\
\ 92 & $ 218.4\leq s<220.8$ & 93 & 64 \ \ \ \  & 55 \ \ \ \  \\
\ 93 & $220.8\leq s<223.2$ & 90  & 61 \ \ \ \  & 51 \ \ \ \  \\
\ 94 & $223.2\leq s<225.6$ & 88 & 58 \ \ \ \  & 50 \ \ \ \  \\
\ 95 & $225.6\leq s<228$    & 86 & 56 \ \ \ \  & 47  \ \ \ \  \\
\ 96 & $228\leq s<230.4$    & 83  & 54\ \ \ \  & 53\ \ \ \ 
\end{tabular}

\end{center}

\newpage

\begin{center}
{\Large TABLE\ 5}
\bigskip

\begin{tabular}{ccc}
$Section$ \ $s$ & $U_{1}(s)$ & $U_{1}(s)$ \\
& $according$ $to$ $(44)$ & $according$ $experiments$ \\
0 & 42 & 41 \\
1 & 45 & 43 \\
2 & 49 & 46 \\
3 & 52 & 48 \\
4 & 56 & 51 \\
5 & 60 & 58 \\
6 & 64 & 65 \\
7 & 68 & 73 \\
8 & 72 & 79 \\
9 & 76 & 83 \\
10 & 80 & 85 \\
11 & 84 & 88 \\
12 & 89 & 91 \\
13 & 90 & 91 \\
14 & 91 & 93 \\
15 & 91 & 94 \\
16 & 92 & 94 \\
17 & 92 & 94 \\
18 & 92 & 94 \\
19 & 92 & 94 \\
20 & 93 & 94 \\
21 & 94 & 94 \\
22 & 94 & 94 \\
23 & 94 & 94 \\
24 & 94 & 93 \\
25 & 95 & 94 \\
26 & 96 & 93 \\
27 & 96 & 92 \\
28 & 96 & 92 \\
29 & 96 & 93 \\
30 & 97 & 95 \\
31 & 97 & 96 \\
32 & 97 & 99 \\
33 & 97 & 100 \\
34 & 97 & 100 \\
35 & 98 & 102 \\
36 & 98 & 102 \\
37 & 98 & 103 \\
38 & 99 & 104 \\
39 & 99 & 104 \\
40 & 99 & 104 \\
41 & 99 & 103 \\
42 & 100 & 102 \\
43 & 100 & 102 \\
44 & 100 & 104 \\
45 & 100 & 98 \\
46 & 100 & 96 \\
47 & 101 & 94 \\
48 & 101 & 91 \\
49 & 101 & 92 \\
50 & 101 & 93 \\
51 & 101 & 93 \\

\end{tabular}

\newpage 
\begin{tabular}{ccc}
$Section$ \ $s$ & $U_{1}(s)$ & $U_{1}(s)$ \\
& $according$ $to$ $(44)$ & $according$ $experiments$ \\
52 & 101 & 95 \\
53 & 101 & 96 \\
54 & 102 & 98 \\
55 & 102 & 99 \\
56 & 102 & 98 \\
57 & 102 & 100 \\
58 & 102 & 103 \\
59 & 103 & 104 \\
60 & 103 & 108 \\
61 & 103 & 109 \\
62 & 103 & 108 \\
63 & 104 & 109 \\
64 & 104 & 108 \\
65 & 105 & 106 \\
66 & 106 & 104 \\
67 & 106 & 103 \\
68 & 106 & 103 \\
69 & 106 & 103 \\
70 & 106 & 104 \\
71 & 107 & 106 \\
72 & 107 & 107 \\
73 & 107 & 106 \\
74 & 107 & 105 \\
75 & 107 & 105 \\
76 & 107 & 104 \\
77 & 107 & 105 \\
78 & 107 & 104 \\
79 & 107 & 105 \\
80 & 107 & 105 \\
81 & 107 & 106 \\
82 & 107 & 106 \\
83 & 107 & 107 \\
84 & 107 & 107 \\
85 & 107 & 106 \\
86 & 107 & 106 \\
87 & 107 & 106 \\
88 & 107 & 105 \\
89 & 107 & 107 \\
90 & 107 & 107 \\
91 & 107 & 108 \\
92 & 107 & 109 \\
93 & 108 & 109 \\
94 & 108 & 110 \\
95 & 108 & 111 \\
96 & 108 & 111 \\
97 & 108 & 111 \\
98 & 109 & 112 \\
99 & 109 & 112 \\
100 & 108 & 110\\
101 & 108 & 110 \\
102 & 109 & 110 \\
103 & 109 & 110 \\
104 & 109 & 111 \\
105 & 109 & 113 \\
\end{tabular}

\newpage
\begin{tabular}{ccc}
$Section$ \ $s$ & $U_{1}(s)$ & $U_{1}(s)$ \\
& $according$ $to$ $(44)$ & $according$ $experiments$ \\
106 & 109 & 114 \\
107 & 110 & 116 \\
108 & 110 & 118 \\
109 & 110 & 118 \\
110 & 110 & 118 \\
111 & 111 & 118 \\
112 & 113 & 117 \\
113 & 114 & 116 \\
114 & 114 & 116 \\
115 & 114 & 117 \\
116 & 114 & 116 \\
117 & 113 & 116 \\
118 & 113 & 115 \\
119 & 114 & 115 \\
120 & 114 & 114 \\
121 & 114 & 114 \\
122 & 113 & 113 \\
123 & 113 & 114 \\
124 & 112 & 114 \\
125 & 112 & 115 \\
126 & 112 & 114 \\
127 & 112 & 115 \\
128 & 112 & 115 \\
129 & 111 & 114 \\
130 & 111 & 114 \\
131 & 110 & 113 \\
132 & 110 & 111 \\
133 & 110 & 110 \\
134 & 110 & 109 \\
135 & 110 & 108 \\
136 & 109 & 108 \\
137 & 109 & 107 \\
138 & 110 & 107 \\
139 & 108 & 107 \\
140 & 108 & 106 \\
141 & 108 & 106 \\
142 & 108 & 106 \\
143 & 107 & 107 \\
144 & 107 & 108 \\
145 & 107 & 108 \\
146 & 106 & 108 \\
147 & 106 & 107 \\
148 & 105 & 106 \\
149 & 105 & 105 \\
150 & 105 & 104 \\
151 & 104 & 103 \\
152 & 104 & 103 \\
153 & 103 & 102 \\
154 & 103 & 101 \\
155 & 102 & 99 \\
156 & 102 & 97 \\
157 & 101 & 97 \\
158 & 100 & 97 \\
159 & 100 & 95 \\
160 & 99 & 95 \\
\end{tabular}

$\bigskip $%
\begin{tabular}{ccc}
$Section$ \ $s$ & $U_{1}(s)$ & $U_{1}(s)$ \\
& $according$ $to$ $(44)$ & $according$ $experiments$ \\
161 & 99 & 94 \\
162 & 98 & 93 \\
163 & 97 & 93 \\
164 & 97 & 93 \\
165 & 96 & 93 \\
166 & 96 & 94 \\
167 & 95 & 94 \\
168 & 95 & 96 \\
169 & 95 & 96 \\
170 & 94 & 95 \\
171 & 94 & 95 \\
172 & 93 & 93 \\
173 & 92 & 93 \\
174 & 92 & 93 \\
175 & 91 & 92 \\
176 & 91 & 94 \\
177 & 90 & 95 \\
178 & 90 & 95 \\
179 & 89 & 95 \\
180 & 88 & 95 \\
181 & 87 & 93 \\
182 & 86 & 93 \\
183 & 86 & 93 \\
184 & 85 & 92 \\
185 & 85 & 92 \\
186 & 84 & 91 \\
187 & 83 & 90 \\
188 & 83 & 91 \\
189 & 82 & 90 \\
190 & 82 & 89 \\
191 & 81 & 88 \\
192 & 81 & 86 \\
193 & 81 & 84 \\
194 & 80 & 83 \\
195 & 80 & 82 \\
196 & 79 & 80 \\
197 & 79 & 79 \\
198 & 79 & 79 \\
199 & 78 & 78 \\
200 & 78 & 77 \\
201 & 77 & 77 \\
202 & 77 & 76 \\
203 & 76 & 75 \\
204 & 76 & 74 \\
205 & 76 & 74 \\
206 & 75 & 73 \\
207 & 75 & 72 \\
208 & 74 & 71 \\
209 & 74 & 70 \\
210 & 73 & 69 \\
\end{tabular}

$\bigskip $%
\begin{tabular}{ccc}
$Section$ \ $s$ & $U_{1}(s)$ & $U_{1}(s)$ \\
& $according$ $to$ $(44)$ & $according$ $experiments$ \\
211 & 72 & 67 \\
212 & 71 & 65 \\
213 & 71 & 64 \\
214 & 70 & 64 \\
215 & 69 & 62 \\
216 & 69 & 61 \\
217 & 68 & 59 \\
218 & 67 & 57 \\
219 & 66 & 55 \\
220 & 65 & 53 \\
221 & 64 & 52 \\
222 & 63 & 51 \\
223 & 61 & 50 \\
224 & 61 & 51 \\
225 & 60 & 50 \\
226 & 59 & 50 \\
227 & 56 & 43 \\
228 & 53 & 55 \\
229 & 51 & 55 \\
230 & 48 & 42
\end{tabular}
\end{center}


\begin{thebibliography}{99}

\bibitem{Patternformation} Forgacs G, Newman SA. 2005. Pattern formation: segmentation, axes and asymmetry. In: Biological physics of the developing embryo. Chapt. 7: 155-187. Cambridge University Press. NY, USA.

\bibitem{Gilbert} Gilbert S. 2003. Developmental Biology, Seventh Edition, Sinauer Associates, Inc Publishers, Massachusetts.

\bibitem{[1]} Wolpert L. 1998. Principles of Development. Oxford University Press, Oxford.

\bibitem{[2]} Kondo S, Asai R. 1995. A reaction-diffusion wave on the skin of the marine angelfish Pomocanthus. Nature 376: 765-768.

\bibitem{[3]} Murray JD. 1988. How the leopard gets its spots. Sci. Amer. 258: 80-87.

\bibitem{[4]} Richardson MK, Hornbruch A, Wolpert L. 1991. Pigment patterns in neural crest chimaeras constituted from quail and guinea fowl embryos. Dev. Biol. 143: 303-319.

\bibitem{[Schoenwolf]} Schoenwolf GC. 2000. Molecular genetic control of axis patterning during early embryogenesis of vertebrates. Ann N Y Acad Sci. 919: 246-60.

\bibitem{[Towers]} Towers M, Mahood R, Yin Y, Tickle C. 2008. Integration of growth and specification in chick wing digit-patterning. Nature 452(7189):882-886.

\bibitem{[Brusco95]} Brusco A, Pecci Saavedra J, Scicolone G, Flores V. 1995. Development of serotonergic innervation of the chick embryo tectum opticum. Int J Dev Neurosci 13(8):835-843.

\bibitem{[Scicolone95]} Scicolone G, Pereyra-Alfonso S, Brusco A, Pecci Saavedra J, Flores V. 1995. Development of the laminated pattern of the chick embryo tectum opticum. Int J Dev Neurosci 13:845-858.

\bibitem{[PereyraAlfonso97]} Pereyra-Alfonso S, Scicolone G, Ferr\'{a}n JL, Pecci Saavedra J, Flores V. 1997.  Developmental pattern of plasminogen activator activity in the chick optic lobe. Int J Dev Neurosci 15:805-812

\bibitem{[PereyraAlfonso98]} Pereyra-Alfonso S, S\'{a}nchez V, Scicolone G, Ferr\'{a}n JL, Flores V. 1998. Cephalo-caudal gradient of plasminogen activator expression in the developing chick optic lobe.  Proc 5th Braz Symp Extracell Matrix SIMEC 98.

\bibitem{[Sanchez02]} S\'{a}nchez V, Ferr\'{a}n JL, Pereyra-Alfonso S, Scicolone G, Rapacioli M, Flores V. 2002. Developmental changes in the spatial pattern of mesencephalic trigeminal nucleus (Mes5) neuron populations in the developing chick optic tectum. J Comp Neurol 448(4):337-348.

\bibitem{[RodGil05]} Rodr\'{i}guez Gil DJ, Vacotto M, Rapacioli M, Scicolone G, Flores V, Fiszer de Plazas S. 2005. Development and localisation of GABA(A) receptor alpha1, alpha2, beta2 and gamma2 subunit mRNA in the chick optic tectum. J Neurosci Res 81(4):469-480.

\bibitem{[Rapacioli00]} Rapacioli M, Gigola S, DÕAttellis C, Ferr‡n JL, Pereyra-Alfonso S, S\'{a}nchez V, Flores V. 2000. Analysis of cell proliferation as stochastic point process. Mathematics and Computers in Modern Science. World Scientific and Engineering Society Press, New York. pp 153-157.

\bibitem{[Mazzeo04a]} Mazzeo J, Rapacioli M, Perfetto J, Fuentes F, Ortalli L, Scicolone G, S\'{a}nchez V, D'Attellis C, Flores V. 2004. Nonlinear analyses of cell proliferation in the central nervous system reveal stochastic and deterministic components. Conf Proc IEEE Eng Med Biol Soc. 2:857-860.

\bibitem{[Mazzeo08]} Mazzeo J, Rapacioli M, Rodriguez Celin A, Duarte S, Flores V. 2008. Multiscale characteristics of cell proliferation in the developing central nervous system of chick embryos. XVI Conference on Nonequilibrium Statistical Mechanics and Nonlinear Physics (XVI Medyfinol 2008) Punta del Este, Uruguay.

\bibitem{[Rapacioli08]} Rapacioli M. 2008. Biolog\'{i}a te\'{o}rica. An\'{a}lisis de la din\'{a}mica de la proliferaci\'{o}n de c\'{e}lulas neuroepiteliales y de la migraci\'{o}n neuronal postmit\'{o}tica durante el desarrollo del sistema nervioso central. Doctoral Thesis. University of Buenos Aires. School of Medicine. Director: Flores V.

\bibitem{segmentation} O'Farrell PH, Stumpff J, Su TT. 2004. Embryonic cleavage cycles: how is a mouse like a fly? Curr Biol. 14(1): 35-45.

\bibitem{[Raff]} Raff MC. 1996. Size control: the regulation of cell numbers in animal development. Cell. 86(2):173-175.

\bibitem{[Bohnsack]} Bohnsack BL, Hirschi KK. 2004. Red light, green light: signals that control endothelial cell proliferation during embryonic vascular development. Cell Cycle. 3(12):1506-1511.

\bibitem{Johnson} Johnson DG, Walker CL. 1999. Cyclins and cell cycle checkpoints. Annu Rev Pharmacol Toxicol. 39:295-312.

\bibitem{[LaVail]} LaVail JH and Cowan WM. 1971. The development of the chick optic tectum. Autoradiographic studies. Brain Res. 28: 421-441.

\bibitem{[Scicolone06]} Scicolone G, Ortalli AL, Alvarez G, Lopez-Costa JJ, Rapacioli M, Ferr\'{a}n JL, S\'{a}nchez V, Flores V. 2006. Developmental pattern of NADPH-diaphorase positive neurons in chick optic tectum is sensitive to changes in visual stimulation. J Comp Neurol 494(6):1007-1030.

\bibitem{[Shigetani]} Shigetani Y, Funahashi JI, Nakamura H. 1997. En-2 regulates the expression of the ligands for Eph type tyrosine kinases in chick embryonic tectum. Neurosci. Res., 27(3): 211-217.

\bibitem{[Luksch]} Luksch H. 2003. Cytoarchitecture of the Avian Optic Tectum: Neuronal Substrate for Cellular Computation.Rev. Neurosci. 14: 85-106.

\bibitem{[Vacotto]} Vacotto M, Rodr\'{i}guez Gil DJ, Mitridate de Novara A, Fiszer de Plazas S. 2003. Differential and irreversible CNS ontogenic reduction in maximal MK-801 binding site number in the NMDA receptor after acute hypoxic hypoxia. Brain Res. 976(2):202-208.

\bibitem{[Fiszer]} Fiszer de Plazas S, Rapacioli M, Gil DJ, Vacotto M, Flores V. 2007. Acute hypoxia differentially affects the gamma-aminobutyric acid type A receptor alpha(1), alpha(2), beta(2), and gamma(2) subunit mRNA levels in the developing chick optic tectum: Stage-dependent sensitivity. J Neurosci Res. 85(14):3135-3144.

\bibitem{[B3]} Perthame B. 2007. Transport equations in biology.
Frontiers in Mathematics. Birkh\"auser Verlag, Basel.

\bibitem{[B1]} Bouchut F. 2004. Nonlinear stability of finite volume methods for hyperbolic conservation laws and well-balanced schemes for sources. Frontiers in Mathematics. Birkh\"auser Verlag, Basel. 

\bibitem{[B2]} LeVeque, RJ. 2002. Finite volume methods for hyperbolic problems. Cambridge Texts in Applied Mathematics. Cambridge University Press, Cambridge.

\bibitem{[Cowan]} Cowan WM, Martin AH, Wenger E. 1968. Mitotic patterns in the optic tectum of the chick during normal development and after early removal of the optic vesicle. Exp Zool. 169(1): 71-92.

\bibitem{[Alvarez]} \'{A}lvarez IS, Mart\'{i}n-Partido G, Rodr\'{i}guez-Gallardo L, Gonz\'{a}lez-Ramos C, Navascu\'{e}s J. 1989. Cell proliferation during early development of the chick embryo otic anlage: quantitative comparison of migratory and nonmigratory regions of the otic epithelium. J Comp Neurol. 290(2): 278-288.

\bibitem{[Kahane]} Kahane N, Kalcheim C. 1998. Identification of early postmitotic cells in distinct embryonic sites and their possible roles in morphogenesis. Cell Tissue Res. 294(2):297-307.

\bibitem{[Freeman]} Freeman M and Gurdon JB. 2002. Regulatory principles of developmental signaling. Annu. Rev. Cell Dev. Biol. 18: 515-539.

\bibitem{[Wolpert]} Wolpert, L. 1998. Principles of Development. Oxford University Press, Oxford.

\bibitem{[Green]} Green J. 2002. Morphogen gradients, positional information and Xenopus: interplay of theory and experiment. Dev Dyn 225: 392-408.

\bibitem{[Rapacioli01]} Rapacioli M, D«Attellis C, DiMiro A, Spraggon T, Ferr\'{a}n JL, Pereyra Alfonso S, S\'{a}nchez V, Scicolone G, Flores V. 2001. Attempts to mathematically define a developmental gradient. The isthmic organizer and the postmitotic neuronal migration dynamics in the developing central nervous system. Mathematics and simulation with biological applications. World Scientific and Engineering Society Press, New York. pp: 137-142.

\bibitem{[Mazzeo04b]} Mazzeo J, Rapacioli M, Fuentes F, Di Guilmi M, Ortalli AL, DÕAttellis C, Flores V. 2004. Space sequences reveal an organized neuroepithelial cell proliferation in the developing central nervous system. WSEAS Transactions on Biology and Biomedicine 4(1): 441-448.

\end{thebibliography}
\end{document}